\documentclass[11pt]{article}
\usepackage{amsmath, upgreek}
\usepackage{amssymb}
\usepackage{amsfonts}
\usepackage{amsmath,tabu}
\usepackage{cite}
\usepackage{color}
\usepackage[final]{pdfpages}
\usepackage{xcolor} 
\usepackage{amscd}
\usepackage{xspace}
\usepackage{verbatim}
\usepackage{graphicx}
\newcommand{\mysection}[1]{
\section{#1}\setcounter{equation}{0}}
\title{\bf A dynamical system approach to the Chandrasekhar-Hamilton-Jacobi equation
}
\author{{\bf Marie-Fran\c{c}oise Bidaut-V\'eron\footnote{\noindent Institut Denis Poisson,
CNRS UMR 7013, Universit\'e de Tours, 37200 Tours, France. E-mail: veronmf@univ-tours.fr}} \\
 {\bf Laurent V\'eron \footnote{\noindent
Institut Denis Poisson,
CNRS UMR 7013,  Universit\'e de Tours, 37200 Tours, France. E-mail: veronl@univ-tours.fr}}\\[2mm]
}
\date{}

\begin{document}
 \maketitle


\newcommand{\txt}[1]{\;\text{ #1 }\;}
\newcommand{\tbf}{\textbf}
\newcommand{\tit}{\textit}
\newcommand{\tsc}{\textsc}
\newcommand{\trm}{\textrm}
\newcommand{\mbf}{\mathbf}
\newcommand{\mrm}{\mathrm}
\newcommand{\bsym}{\boldsymbol}
\newcommand{\scs}{\scriptstyle}
\newcommand{\sss}{\scriptscriptstyle}
\newcommand{\txts}{\textstyle}
\newcommand{\dsps}{\displaystyle}
\newcommand{\fnz}{\footnotesize}
\newcommand{\scz}{\scriptsize}
\newcommand{\be}{\begin{equation}}
\newcommand{\bel}[1]{\begin{equation}\label{#1}}
\newcommand{\ee}{\end{equation}}
\newcommand{\eqnl}[2]{\begin{equation}\label{#1}{#2}\end{equation}}
\newcommand{\barr}{\begin{eqnarray}}
\newcommand{\earr}{\end{eqnarray}}
\newcommand{\bars}{\begin{eqnarray*}}
\newcommand{\ears}{\end{eqnarray*}}
\newcommand{\nnu}{\nonumber \\}
\newtheorem{subn}{\name}
\renewcommand{\thesubn}{}
\newcommand{\bsn}[1]{\def\name{#1}\begin{subn}}
\newcommand{\esn}{\end{subn}}
\newtheorem{sub}{\name}[section]
\newcommand{\dn}[1]{\def\name{#1}}   
\newcommand{\bs}{\begin{sub}}
\newcommand{\es}{\end{sub}}
\newcommand{\bsl}[1]{\begin{sub}\label{#1}}
\newcommand{\bth}[1]{\def\name{Theorem}
\begin{sub}\label{t:#1}}
\newcommand{\blemma}[1]{\def\name{Lemma}
\begin{sub}\label{l:#1}}
\newcommand{\bcor}[1]{\def\name{Corollary}
\begin{sub}\label{c:#1}}
\newcommand{\bdef}[1]{\def\name{Definition}
\begin{sub}\label{d:#1}}
\newcommand{\bprop}[1]{\def\name{Proposition}
\begin{sub}\label{p:#1}}

\newcommand{\R}{\eqref}
\newcommand{\rth}[1]{Theorem~\ref{t:#1}}
\newcommand{\rlemma}[1]{Lemma~\ref{l:#1}}
\newcommand{\rcor}[1]{Corollary~\ref{c:#1}}
\newcommand{\rdef}[1]{Definition~\ref{d:#1}}
\newcommand{\rprop}[1]{Proposition~\ref{p:#1}}
\newcommand{\BA}{\begin{array}}
\newcommand{\EA}{\end{array}}
\newcommand{\BAN}{\renewcommand{\arraystretch}{1.2}
\setlength{\arraycolsep}{2pt}\begin{array}}
\newcommand{\BAV}[2]{\renewcommand{\arraystretch}{#1}
\setlength{\arraycolsep}{#2}\begin{array}}
\newcommand{\BSA}{\begin{subarray}}
\newcommand{\ESA}{\end{subarray}}
\newcommand{\BAL}{\begin{aligned}}
\newcommand{\EAL}{\end{aligned}}
\newcommand{\BALG}{\begin{alignat}}
\newcommand{\EALG}{\end{alignat}}
\newcommand{\BALGN}{\begin{alignat*}}
\newcommand{\EALGN}{\end{alignat*}}
\newcommand{\note}[1]{\textit{#1.}\hspace{2mm}}
\newcommand{\Proof}{\note{Proof}}
\newcommand{\qeda}{\hspace{10mm}\hfill $\square$}
\newcommand{\qed}{\\
${}$ \hfill $\square$}
\newcommand{\Remark}{\note{Remark}}
\newcommand{\modin}{$\,$\\[-4mm] \indent}
\newcommand{\forevery}{\quad \forall}
\newcommand{\set}[1]{\{#1\}}
\newcommand{\setdef}[2]{\{\,#1:\,#2\,\}}
\newcommand{\setm}[2]{\{\,#1\mid #2\,\}}
\newcommand{\mt}{\mapsto}
\newcommand{\lra}{\longrightarrow}
\newcommand{\lla}{\longleftarrow}
\newcommand{\llra}{\longleftrightarrow}
\newcommand{\Lra}{\Longrightarrow}
\newcommand{\Lla}{\Longleftarrow}
\newcommand{\Llra}{\Longleftrightarrow}
\newcommand{\warrow}{\rightharpoonup}
\newcommand{
\paran}[1]{\left (#1 \right )}
\newcommand{\sqbr}[1]{\left [#1 \right ]}
\newcommand{\curlybr}[1]{\left \{#1 \right \}}
\newcommand{\abs}[1]{\left |#1\right |}
\newcommand{\norm}[1]{\left \|#1\right \|}
\newcommand{
\paranb}[1]{\big (#1 \big )}
\newcommand{\lsqbrb}[1]{\big [#1 \big ]}
\newcommand{\lcurlybrb}[1]{\big \{#1 \big \}}
\newcommand{\absb}[1]{\big |#1\big |}
\newcommand{\normb}[1]{\big \|#1\big \|}
\newcommand{
\paranB}[1]{\Big (#1 \Big )}
\newcommand{\absB}[1]{\Big |#1\Big |}
\newcommand{\normB}[1]{\Big \|#1\Big \|}
\newcommand{\produal}[1]{\langle #1 \rangle}

\newcommand{\thkl}{\rule[-.5mm]{.3mm}{3mm}}
\newcommand{\thknorm}[1]{\thkl #1 \thkl\,}
\newcommand{\trinorm}[1]{|\!|\!| #1 |\!|\!|\,}
\newcommand{\bang}[1]{\langle #1 \rangle}
\def\angb<#1>{\langle #1 \rangle}
\newcommand{\vstrut}[1]{\rule{0mm}{#1}}
\newcommand{\rec}[1]{\frac{1}{#1}}
\newcommand{\opname}[1]{\mbox{\rm #1}\,}
\newcommand{\supp}{\opname{supp}}
\newcommand{\dist}{\opname{dist}}
\newcommand{\myfrac}[2]{{\displaystyle \frac{#1}{#2} }}
\newcommand{\myint}[2]{{\displaystyle \int_{#1}^{#2}}}
\newcommand{\mysum}[2]{{\displaystyle \sum_{#1}^{#2}}}
\newcommand {\dint}{{\displaystyle \myint\!\!\myint}}
\newcommand{\q}{\quad}
\newcommand{\qq}{\qquad}
\newcommand{\hsp}[1]{\hspace{#1mm}}
\newcommand{\vsp}[1]{\vspace{#1mm}}
\newcommand{\ity}{\infty}
\newcommand{\prt}{\partial}
\newcommand{\sms}{\setminus}
\newcommand{\ems}{\emptyset}
\newcommand{\ti}{\times}
\newcommand{\pr}{^\prime}
\newcommand{\ppr}{^{\prime\prime}}
\newcommand{\tl}{\tilde}
\newcommand{\sbs}{\subset}
\newcommand{\sbeq}{\subseteq}
\newcommand{\nind}{\noindent}
\newcommand{\ind}{\indent}
\newcommand{\ovl}{\overline}
\newcommand{\unl}{\underline}
\newcommand{\nin}{\not\in}
\newcommand{\pfrac}[2]{\genfrac{(}{)}{}{}{#1}{#2}}

\def\ga{\alpha}     \def\gb{\beta}       \def\gg{\gamma}
\def\gc{\chi}       \def\gd{\delta}      \def\ge{\epsilon}
\def\gth{\theta}                         \def\vge{\varepsilon}
\def\gf{\phi}       \def\vgf{\varphi}    \def\gh{\eta}
\def\gi{\iota}      \def\gk{\kappa}      \def\gl{\lambda}
\def\gm{\mu}        \def\gn{\nu}         \def\gp{\pi}
\def\vgp{\varpi}    \def\gr{\rho}        \def\vgr{\varrho}
\def\gs{\sigma}     \def\vgs{\varsigma}  \def\gt{\tau}
\def\gu{\upsilon}   \def\gv{\vartheta}   \def\gw{\omega}
\def\gx{\xi}        \def\gy{\psi}        \def\gz{\zeta}
\def\Gg{\Gamma}     \def\Gd{\Delta}      \def\Gf{\Phi}
\def\Gth{\Theta}
\def\Gl{\Lambda}    \def\Gs{\Sigma}      \def\Gp{\Pi}
\def\Gw{\Omega}     \def\Gx{\Xi}         \def\Gy{\Psi}

\def\CS{{\mathcal S}}   \def\CM{{\mathcal M}}   \def\CN{{\mathcal N}}
\def\CR{{\mathcal R}}   \def\CO{{\mathcal O}}   \def\CP{{\mathcal P}}
\def\CA{{\mathcal A}}   \def\CB{{\mathcal B}}   \def\CC{{\mathcal C}}
\def\CD{{\mathcal D}}   \def\CE{{\mathcal E}}   \def\CF{{\mathcal F}}
\def\CG{{\mathcal G}}   \def\CH{{\mathcal H}}   \def\CI{{\mathcal I}}
\def\CJ{{\mathcal J}}   \def\CK{{\mathcal K}}   \def\CL{{\mathcal L}}
\def\CT{{\mathcal T}}   \def\CU{{\mathcal U}}   \def\CV{{\mathcal V}}
\def\CZ{{\mathcal Z}}   \def\CX{{\mathcal X}}   \def\CY{{\mathcal Y}}
\def\CW{{\mathcal W}} \def\CQ{{\mathcal Q}}
\def\BBA {\mathbb A}   \def\BBb {\mathbb B}    \def\BBC {\mathbb C}
\def\BBD {\mathbb D}   \def\BBE {\mathbb E}    \def\BBF {\mathbb F}
\def\BBG {\mathbb G}   \def\BBH {\mathbb H}    \def\BBI {\mathbb I}
\def\BBJ {\mathbb J}   \def\BBK {\mathbb K}    \def\BBL {\mathbb L}
\def\BBM {\mathbb M}   \def\BBN {\mathbb N}    \def\BBO {\mathbb O}
\def\BBP {\mathbb P}   \def\BBR {\mathbb R}    \def\BBS {\mathbb S}
\def\BBT {\mathbb T}   \def\BBU {\mathbb U}    \def\BBV {\mathbb V}
\def\BBW {\mathbb W}   \def\BBX {\mathbb X}    \def\BBY {\mathbb Y}
\def\BBZ {\mathbb Z}

\def\GTA {\mathfrak A}   \def\GTB {\mathfrak B}    \def\GTC {\mathfrak C}
\def\GTD {\mathfrak D}   \def\GTE {\mathfrak E}    \def\GTF {\mathfrak F}
\def\GTG {\mathfrak G}   \def\GTH {\mathfrak H}    \def\GTI {\mathfrak I}
\def\GTJ {\mathfrak J}   \def\GTK {\mathfrak K}    \def\GTL {\mathfrak L}
\def\GTM {\mathfrak M}   \def\GTN {\mathfrak N}    \def\GTO {\mathfrak O}
\def\GTP {\mathfrak P}   \def\GTR {\mathfrak R}    \def\GTS {\mathfrak S}
\def\GTT {\mathfrak T}   \def\GTU {\mathfrak U}    \def\GTV {\mathfrak V}
\def\GTW {\mathfrak W}   \def\GTX {\mathfrak X}    \def\GTY {\mathfrak Y}
\def\GTZ {\mathfrak Z}   \def\GTQ {\mathfrak Q}

\font\Sym= msam10 
\def\SYM#1{\hbox{\Sym #1}}
\newcommand{\bdw}{\prt\Gw\xspace}
\date{}
\maketitle\smallskip

\begin{center}
\it To Ha\"im Brezis, outstanding mathematician, great mind and sincere friend
\end{center}\medskip

\noindent{\small {\bf Abstract} We study the local properties of positive solutions of the equation $-\Gd u=e^u-M\abs{\nabla u}^q$ in a punctured domain $\Gw\setminus\{0\}$ of $\BBR^N$ in the range of parameters $q>1$ and $M> 0$. We prove a series of a priori estimates near a singular point. In the case of radial solutions we use various techniques inherited from the dynamical systems 
theory to obtain the precise behaviour of singular solutions. We prove also the existence of singular solutions with these precise behaviours.
}\medskip

\tableofcontents
\vspace{1mm}
\hspace{.05in}
\medskip

\mysection{Introduction}
The aim of this paper is to study properties of {\it radial} solutions of 
\bel{I-1}
\CL_{m,p,q}u:=-\Gd u +M|\nabla u|^q-e^u=0\quad\text{in }\Gw,
\ee
where $M$ is a real positive number, $q> 1$ and $\Gw$ is either a punctured domain if we are interested in isolated singularities, or an exterior domain if we study the asymptotic behaviour of solutions. 
\smallskip 
Equation (\ref{I-1}) belongs to the family of equations  of Diffusion-Reaction-Absorption type
\bel{I-5-1}
-\Gd u \pm f(|\nabla u|)=g(u),
\ee
where $f$ is a power function and $g$ a power or an exponential function.\smallskip

For these equations, there are several deep problems:\\
- Obtention of a priori estimates near a singularity or at infinity for equations in an exterior domain.\\
-  Description of singular solutions or behaviour at infinity of solutions.\\
- Existence of global regular or singular solutions. \smallskip

Equations with absorption correspond to the cases where $g(u)=-e^u$ or $g(u)=-u^p$. The equation 
$$
-\Gd u +M |\nabla u|^q+u^p=0,
$$
have been studied in \cite{Ng1} when $M>0$ and in \cite{BiGaVe4} when $M<0$. The equation 
$$
-\Gd u +M |\nabla u|^q+e^{u}=0,
$$
with $M<0$ is considered in \cite{LiVe}, particulary in the most relevant dimension $N=2$ since in higher dimension the isolated singularities are removable.\smallskip

Next we focus on the equation where $g(u)$ acts as a source term. Equation with a power $u^p$ with $p>0$
\bel{I-5-2}
-\Gd u+M|\nabla u|^q- u^p=0,
\ee
has been first considered with $M=1$ in \cite{ChWe1} and \cite{ChWe2} in the radial case where two critical values appear 
$p=\frac{N+2}{N-2}$ and $q=\frac{2p}{p+1}$. This study was developed in \cite{SeZo} and \cite{SeZo2}, putting into light the question of existence of global solutions, still not completely solved. The general problem with $M\in\BBR$ has been studied in the non-radial case in \cite{PoQuSo} when $p<\frac{N+2}{N-2}$ and $q<\frac{2p}{p+1}$, then in \cite{BiGaVe2} when $q>\frac{2p}{p+1}$. The critical case $q=\frac {2p}{p+1}$ with $p>1$ has been exhaustively solved in the radial case in \cite{BiGaVe3}.
\smallskip

In the present article we assume that $g(u)$ is an exponential source term and we consider equation $(\ref{I-1})$ with $M>0$. To study the properties of the solutions of ${\ref{I-1}}$ we associate three underlying equations:\\
the {\it Emden-Chandrashekar equation}
\bel{I-E-C}
-\Gd u-e^{u}=0,
\ee
the {\it viscous Hamilton-Jacobi equation}
\bel{I-HJ}
-\Gd u+M|\nabla u|^q=0,
\ee
and the {\it eikonal equation}
\bel{I-ei}
M|\nabla u|^q-e^{u}=0.
\ee
Each of this three  equations have been already much studied. \\
- The first exhaustive study of radial solutions of the Emden-Chandrashekar equation is due to Chandrashekar \cite{Cha}. Non-radial solutions which are rather simple to study 
in the 2-dimensional case are described in \cite{RiVe}. In the 3-dimensional case the level of difficulties is not at all comparable and this is due to the fact that the reaction term is much larger than the classical Sobolev exponent $p=\frac{N+2}{N-2}$ which turns out to be equal to $5$ since $N=3$. It is proved in \cite{BiVe0} that no uniform estimate near a singularity can hold, but if a solution $u$ 
of $(\ref{I-E-C})$ in $B_1\setminus\{0\}$ (resp. $B^c_1$) satisfies $|x|^2e^{u}\in L^\infty(B_1)$ (resp.  $|x|^2e^{u}\in L^\infty(B^c_1)$) then its behaviour near $x=0$ (resp. when $|x|\to\infty$) can be completely described. Besides the upper estimate another severe problem of convergence occurs and is overcame by the introduction of deep results from complex geometry.\\
- The viscous Hamilton-Jacobi equation is simpler to study. In the radial case it reduces to a explicitely integrable equation
$$-u_{rr}-\frac{N-1}{r}u_r+M|u_r|^q=0.
$$
This equation admits two critical exponents $q=q_c=\frac{N}{N-1}$ and $q=2$. By using Bernstein method, it is proved in \cite{La-Li} that if $q>1$ any not necessarily radial solution $u$ of $(\ref{I-HJ})$ in $B_1\setminus\{0\}$ satisfies 
$$|\nabla u(x)|\leq c|x|^{-\frac{1}{q-1}}\quad\text{in }B_{\frac 12}\setminus\{0\},
$$
and that its behaviour near zero is similar as the one of the radial (and explicit) solutions, in particular 
if $q\geq q_c$ any nonnegative solution remains bounded, while if $1<q<q_c$ and if it is singular either it satisfies $u(x)\sim c|x|^{2-N})$ for some $c\neq 0$, or  $u(x)\sim \xi_M|x|^{\frac{q-2}{q-1}}$, where $\xi_M$ is a specific positive constant expressed in $(\ref{IV-6})$.\\
- The eikonal equation is also explicitely integrable in the radial case and always it reduces to an equation of the form
$$|\nabla v|=1
$$
with $v=qM^{\frac 1q}e^{-\frac uq}$ in general. The uniqueness of viscosity solutions of the eikonal equation in $\BBR^N\setminus\{0\}$ is proved in \cite{Ca-Cr}.\smallskip 

In \cite{BiVe2} we present a general study of non-necessarily radial solutions of $(\ref{I-1})$ and emphasis on their singular solutions if they are defined in $B_1\setminus\{0\}$, and their asymptotic behaviour if they defined in 
$B^c_1$. This study involves all the tools used in \cite{BiVe0}, often under a more refined manner, and $q=2$ appears as a fundamental critical value. In particular if $1<q<2$ we have to impose the a priori bound
$|x|^2e^u\leq C$ in $B_1\setminus\{0\}$ in order to obtain the description of singular solutions, while if $q>2$ no a priori bound on a solution is needed. \smallskip

In the present article we consider functions satisfying 
\bel{I-7}
-u_{rr}-\frac{N-1}{r}u_r+M|u_r|^q-e^u=0,
\ee
 with $q>1$ and $M>0$, either in $(0,1]$ or in $[1,\infty)$. This assumption of radiality has the advantage of authorising the use of finite dimensional dynamical systems theory, a theory which allows us to go much deeper in the local study of the solutions. For $q\neq 2$ there is no  invariance for $(\ref{I-7})$ 
by a scaling transformation (the case $q=2$ which reduces to a Lane-Emden type equation will not be considered). Consequently it is not possible to reduces this equation to a second order autonomous equation of any type as it is possible for $(\ref{I-E-C})$, $(\ref{I-HJ})$ and $(\ref{I-ei})$, and also for $(\ref{I-5-2})$ when $q=\frac{2p}{p+1}$. 
Instead we introduce various systems of order 3 to perform our analysis of singular and asymptotic behaviours: we  give a precise description of singular solutions (or there asymptotic behaviour), and we prove the existence of solutions 
with all the possible behaviours that we have put into light. These existence results were out of reach for non-necessarily radial solutions studied in  \cite{BiVe2}.\medskip

An important observation, dealing  with the globality of the radial solutions is the following statement proved in \rprop{glob} and \rth{globsub}.\medskip

\nind{\bf Theorem A} {\it  Let $q>1$.\smallskip

\nind (1) If  $N\geq 2$, the maximal interval of definition of a radial solution $r\mapsto u(r)$ of $(\ref{I-7})$ defined  in a right neighborhood of  $r=0$ and nondecreasing there  is $(0,\infty)$. Furthermore $u_r<0$ and $u$ tends to $-\infty$ at $\infty$.\smallskip

\nind (2) If $1<q<2$ and $N\geq 2$, the maximal interval of definition of a solution $r\mapsto u(r)$ of $(\ref{I-7})$  is $(0,\infty)$. If $1<q<2$ and $N=1$ any maximal solution of $(\ref{I-7})$ is  defined on $\BBR$ and symmetric with respect to some value $a$.}\smallskip

 We obtain also a full description of the isolated singularities of solutions of $(\ref{I-7})$.\medskip

 \nind {\bf Theorem B} {\it Let $1<q<2$ and $N\geq 3$, \smallskip
 
 \nind (1) either 
$\dsps \lim_{r\to 0}\left(u(r)-2\ln\frac 1r\right)=\ln 2(N-2),$\\
\nind (2) or there exists $u_0\in\BBR$ such that $\dsps\lim_{r\to 0}u(r)=u_0$ and $\dsps\lim_{r\to 0}u_r(r)=0$ (such a $u$ is called a regular solution),\\
\nind (3) or $\dsps\lim_{r\to 0}u(r)=-\infty$ and:\smallskip

\nind - if  $1<q<\frac{N}{N-1}$, then $\dsps\lim_{r\to 0}r^{N-2}u(r)=\gg<0$,\\
- or  if $q=\frac{N}{N-1}$, then $\dsps \lim_{r\to 0}r^{N-2}|\ln r|^{N-1}u(r)=-\frac{1}{N-2}\left(\frac{N-1}{M}\right)^{N-1}$,\\
- or  if $\frac{N}{N-1}<q<2$, then $\dsps \lim_{r\to 0}r^{\frac{2-q}{q-1}}u(r)=-\xi_M:=-\frac{q-1}{2-q}\left(\frac{(N-1)q-N}{M(q-1)}\right)^{\frac{1}{q-1}}$.
}\\

This result shows that the singular behaviour of solutions is governed either by the Emden-Chandrashekar equation or by the viscous Hamilton-Jacobi equation. Existence of solutions with the above behaviour is obtained in \rth{nega} and \rth{dirac} by fixed point argument or by the study near the equilibrium of the associated dynamical system.\smallskip

When $q>2$ the situation changes completely
and we prove that the singularities of radial solutions are either governed by the eikonal equation or by the Hamilton-Jacobi equation.  \medskip
 
 \nind  \nind {\bf Theorem C} {\it  Let $q>2$ and $N\geq 2$, \smallskip
 
 \nind (1) either 
$\dsps \lim_{r\to 0}r^qe^{u(r)}=Mq^q$ and $\dsps \lim_{r\to 0}ru_r(r)=-q$,\\
\nind (2) or there exists $u_0\in\BBR$ such that $\dsps\lim_{r\to 0}u(r)=u_0$ and $\dsps\lim_{r\to 0}u_r(r)=0$, in that case $u$ is a regular solution,\\
\nind (3) or there exists $u_0\in\BBR$ such that $\dsps\lim_{r\to 0}u(r)=u_0$ and there  holds
$u(r)=u_0+C_{M,N}r^{\frac{q-2}{q-1}}(1+o(1))$ where $C_{M,N}$ is explicited in Theorem E. In such a case the singularity is only on the derivative.
}\\

We also obtain the existence of singular solutions with the prescribed behaviour given above. Concerning solutions of eikonal type 
we prove existence in \rth{uniq-ei} and \rth{uniq-ei-2} by a very delicate method based upon inverse function arguments.\smallskip

\nind {\bf Theorem D} {\it Let $q>2$. \smallskip

\nind (1) If $N=1$, there exists one and only one solution $u^*$ of $(\ref{I-7})$ on $(0,\infty)$ such that 
\bel{e33}\lim_{r\to 0}r^qe^{u^*(r)}=Mq^q\quad\text{and }\, \lim_{r\to 0}ru^*_r(r)=-q.
\ee
Furthermore the function $u^*$ is the increasing limit when $n\to\infty$ of the regular solutions $u_n$ (i.e. $u_n(0)=n$ and $u_{n\,r}=0$)}.\smallskip

\nind {\it  (2) If $N\geq 2$ there exists at least one solution $u^*$ of $(\ref{I-7})$ on $(0,\infty)$ satisfying $(\ref{e33})$.
}\\

Concerning solutions of Hamilton-Jacobi type, we prove,\smallskip

\nind {\bf Theorem E} {\it Let $q>2$ and $u_0 \in \BBR$ arbitrary. \smallskip

\nind (1) If $N=1$ there exists at least one solution of $(\ref{I-7})$ on $(0,\infty)$ satisfying 
$$u(r)=u_0-\frac{q-1}{q-2}\left(\frac{1}{M(q-1)}\right)^{\frac1{q-1}}r^{\frac{q-2}{q-1}}(1+o(1))\quad\text{as }r\to 0.
$$
Furthermore the function $u$ is decreasing on $(0,\infty)$}.\smallskip

\nind {\it  (2) If $N\geq 2$ there exists at least one solution of of $(\ref{I-7})$ on $(0,\infty)$ satisfying 
$$u(r)=u_0+\frac{q-1}{q-2}\left(\frac{(N-1)q-N}{M(q-1)}\right)^{\frac1{q-1}}r^{\frac{q-2}{q-1}}(1+o(1))\quad\text{as }r\to 0.
$$
Furthermore the function $u$ is increasing on $(0,\infty)$.
}\\

 This result is proved in \rth{q>2-beta} and \rth {q>2-beta*} by methods coming from the analysis of the stable and unstable manifolds associated to the  stationary points of the relevant dynamical system.\smallskip

The description of the asymptotic behaviour of radial solutions of $(\ref{I-1})$ in an exterior domain exchanges the ranges $1<q<2$ and $q>2$. Note also that only one type of 
behaviour is possible. The following statements are proved in \rth{inftyq<2} and \rth{inftyq>2}.\smallskip

\nind \nind {\bf Theorem F} {\it  (1) Let $1<q<2$ and $N\geq 3$. If $u$ is a radial solution of $(\ref{I-1})$ in $\overline B_{r_0^c}$ it satisfies $\dsps\lim_{r\to\infty}r^qe^{u(r)}=Mq^q$ and $\dsps \lim_{r\to 0}ru_r(r)=-q$.}\smallskip

\nind{\it (2) Let $q>2$ and $N\geq 3$. If $u$ is a radial solution of $(\ref{I-1})$ in $\overline B_{r_0^c}$ it satisfies $\dsps\lim_{r\to\infty}r^2e^{u(r)}=2(N-2)$.}\\

We end the article with an appendix which shows that when $q>2$, all the radial solutions $u$ of $(\ref{I-1})$ which  satisfy
$\dsps \lim_{r\to 0}r^qe^{u(r)}=Mq^q$ and $\dsps \lim_{r\to 0}ru_r(r)=-q$, which means that they are of eikonal type, have the same expansion at any order near $0$. This is a 
clue that uniqueness could hold in any dimension as it does hold when $N=1$.\\

We leave as an open problem the study of 
\bel{BM1}
-\Gd u+M|\nabla u|^q=V(x)e^u.
\ee
for many types of potential $V(x)$. Deep results in the case $N=2$, $M=0$ have been obtained by Brezis and Merle \cite{BrMe}.

\mysection{Estimates}
\subsection{Estimates of radial supersolutions}

\bth{sup-1} Let $N\geq 1$ and $q>1$.\smallskip

\nind (i) If $u\in C(\overline {B_{r_0}}\setminus\{0\})$ is a radial supersolution of $(\ref{I-1})$ in $B_{r_0}\setminus\{0\}$ it satisfies
\bel{III-2}
e^{u(r)}\leq Cr^{-\max\{2,q\}}\quad\text{in }B_{\frac{r_0}{2}}\setminus\{0\},
\ee
where $C=C(N,M,q,u)>0)$ in the general case and $C=C(N,M,q)>0)$ if $r\mapsto u(r)$ is non-increasing.\smallskip
 
\nind (ii) If $u\in C(\overline {B^c_{r_0}})$ is a supersolution of $(\ref{I-1})$ in $B^c_{r_0}$ it satisfies
\bel{III-3}
e^{u(r)}\leq Cr^{-\min\{2,q\}}\quad\text{in }B^c_{2r_0},
\ee
where $C=C(N,M,q,u)>0$ is as in case (i).
\es
\Proof If $u$ is a radial supersolution it is clear that it has at most one local maximum. Indeed at each local extremal point $\tilde r$ there holds $-u_{rr}(\tilde r)\geq e^{u(\tilde r)}>0$. Hence $u_r$ keeps a constant sign near $0$.
We set $w=e^u$, then $w\geq 0$ satisfies.

\bel{III-5}
-\Gd w+\frac{|\nabla w|^2}{w}+M\frac{|\nabla w|^q}{w^{q-1}}-w^2
=-w_{rr}-\frac{N-1}{r}w_r+\frac{w_r^2}{w}+M\frac{|w_r|^q}{w^{q-1}}-w^2\geq 0.
\ee

\nind (i) We first assume that  the function $w$ is nondecreasing on $(0,r_1]$ for some $r_1\in (0,r_0]$, the estimate $(\ref{III-2})$ holds with $C$ depending on $u$.\\
Next we assume that $m$ is nonincreasing on $(0,r_1]$, then it stays decreasing on the whole interval $(0,r_0]$. For $\ge\in (0,\frac 12)$ let $\phi_\ge$ be a $C^\infty(\BBR_+)$ nonnegative function vanishing on $[0,1-\ge]\cup [1+\ge,\infty)$,
with value $1$ on $[1-\frac\ge2,1+\frac\ge 2]$, such that $\ge|\phi_{\ge\,r}|$ and $\ge^2|\phi_{\ge\,rr}|$ are bounded on $[1-\ge,1-\frac\ge 2]\cup[1+\frac\ge 2,1+\ge]$. If $0<R<\frac{r_0}{2}$ we define $v$ by
$$v(x)=w(x)-w(R)\phi_\ge\left(\frac{|x|}{R}\right)
$$
Clearly  $v(R)=0$. There exists $\tilde r_{R,\ge}\in [1-\ge)R,1+\ge)R]$ where 
$v$ achieves a nonpositive minimum, hence $w(\tilde r_{R,\ge})\leq w(R)$, $v_r(\tilde r_{R,\ge})=0$ and $\Gd v(\tilde r_{R,\ge})\geq 0$. There holds 
$w_r(\tilde r_{R,\ge})=w(R)\phi_{\ge\,r}(\frac {r_{R,\ge}}{R})$ and $\Gd w(\tilde r_{R,\ge})\geq w(R)\Gd \phi_\ge(\frac{|x_{R,\ge}|}{R})$. Therefore
\bel{III-6}\BA{lll}\dsps
w^2(\tilde r_{R,\ge})\leq -\Gd w(\tilde r_{R,\ge})+\frac{|\nabla w(\tilde r_{R,\ge})|^2}{w(\tilde r_{R,\ge})} +M\frac{|\nabla w(\tilde r_{R,\ge})|^q}{w^{q-1}(\tilde r_{R,\ge})} \\[4mm]
\phantom{w^2(\tilde r_{R,\ge})}
\dsps \leq -w(R)\Gd \phi_\ge(\frac{r_{R,\ge}}{R})+w^2(R)\frac{|\nabla \phi_\ge(\frac{r_{R,\ge}}{R})|^2}{w(\tilde r_{R,\ge})} +
w^q(R)\frac{|\nabla \phi_\ge(\frac{r_{R,\ge}}{R})|^q}{w^{q-1}(\tilde r_{R,\ge})}\\[4mm]
\phantom{w^2(\tilde r_{R,\ge})}\dsps
\leq C\left(\frac{w(R)}{\ge^2R^2}+\frac{w^2(R)}{\ge^2R^2w(\tilde r_{R,\ge})}+\frac{w^q(R)}{\ge^qR^qw^{q-1}(\tilde r_{R,\ge})}\right).
\EA\ee
where $C=C(N,q,M)>0$. \\
Now, if $q>2$, we multiply by $w^{q-1}(\tilde r_{R,\ge})$ and obtain the estimate
$$\BA{lll}\dsps
w^{q+1}(\tilde r_{R,\ge})\leq C\left(\frac{w(R)}{\ge^2 R^2}w^{q-1}(R)+\frac{w^2(R)}{\ge^2 R^2}w^{q-2}(R)+\frac{w^q(R)}{\ge^q R^q}\right)\leq C'\frac{w^q(R)}{\ge^q R^q}.
\EA$$
Because $\dsps w^{q+1}(\tilde r_{R,\ge})\geq \min_{(1-\ge)R\leq r\leq (1+\ge)R}w^{q+1}(r)\geq w^{q+1}((1+\ge)R)$ we obtain
$$w^{q+1}((1+\ge)R)\leq C\frac{w^{\frac q{q+1}}(R)}{\ge^{\frac{q}{q+1}}R^{\frac q{q+1}}}.
$$
We apply the bootstrap method of \cite[Lemma 2.1]{BiVe1} with $\Phi(\gr)=\gr^{-\frac q{q+1}}$  and $d=-h=\frac q{q+1}$ and we conclude that 
$$w(R)\leq C(\Phi(R))^{\frac{1}{1-d}}=CR^{-q}.
$$
If $1<q<2$ we have from $(\ref{III-6})$, 
$$\BA{lll}\dsps w^3(\tilde r_{R,\ge})\leq C\left(\frac{w(R)}{\ge^2R^2}w(\tilde r_{R,\ge})+\frac{w^2(R)}{\ge^2R^2}+\frac{w^q(R)}{\ge^qR^qw^{2-q}}(\tilde r_{R,\ge})\right)\\[4mm]
\phantom{w^3(\tilde r_{R,\ge})}\dsps
\leq C\left(\frac{w(R)}{\ge^2R^2}w(R)+\frac{w^2(R)}{\ge^2R^2}+\frac{w^q(R)}{\ge^qR^qw^{2-q}(r})\right)
\\[4mm]
\phantom{w^3(\tilde r_{R,\ge})}\dsps
\leq C'\left(\frac{w^2(R)}{\ge^2R^2}+\frac{w^2(R)}{\ge^qR^2}\right)\leq C''\frac{w^2(R)}{\ge^2R^2}.
\EA$$
Then 
$$w((1+\ge)R)\leq c\frac{w^{\frac 23}(R)}{\ge^{\frac 23}R^\frac 23}.
$$
Now applying \cite[Lemma 2.1]{BiVe1} with $\Phi(\gr)=\gr^{-\frac 23}$, $d=-h=\frac 23$, we obtain
$$w(R)\leq C(\Phi(R))^{\frac{1}{1-d}}=CR^{-2}.
$$
Note that in the two cases, the upper estimate is independent of $u$.\smallskip

\nind (ii) Assume that $u$ satisfies  $(\ref{III-2})$ in $B^c_{r_0}$. Here also the function $w$ is monotone on $[r_1,\infty)$ for some $r_1\geq r_0$. Let $0<\ge\leq\frac 12$ and $R>r_1+1$. If $w$ is nondecreasing on $[R(1-\ge),R(1+\ge)]$, then 
$w(R)\geq w(\tilde r_{R,\ge})\geq w(R(1-\ge))$ and therefore
$$w(R(1-\ge))\leq w(\tilde r_{R,\ge})\leq C\left(\frac{w^q(R)}{\ge^2R^2}+\frac{w^q(R)}{\ge^qR^q}\right)^{\frac 1{q+1}}
\leq\left\{\BA{lll}\dsps2C\frac{w^{\frac{q}{q+1}}}{\ge^{\frac{2}{q+1}}R^{\frac{2}{q+1}}}\quad\text{if }q>2\\[4mm]
\dsps2C\frac{w^{\frac{q}{q+1}}}{\ge^{\frac{q}{q+1}}R^{\frac{q}{q+1}}}\quad\text{if }1<q<2.
\EA\right.
$$
From the bootstrap argument,
$$w(R)\leq\left\{\BA{lll} C_1R^{-2}\quad&\text{if }q>2\\[2mm]
2CR^{-q}\quad&\text{if }1<q<2,
`\EA\right.
$$
which contradicts the fact that $w$ is nondecreasing. Hence $w$ is nonincreasing and we obtain the estimate
$$w(R(1+\ge))\leq\left\{\BA{lll}\dsps2C\frac{w^{\frac{q}{q+1}}}{\ge^{\frac{2}{q+1}}R^{\frac{2}{q+1}}}\quad\text{if }q>2\\[4mm]
\dsps2C\frac{w^{\frac{q}{q+1}}}{\ge^{\frac{q}{q+1}}R^{\frac{q}{q+1}}}\quad\text{if }1<q<2.
\EA\right.$$
Therefore we obtain $(\ref{III-2})$ as before.\qeda\medskip

\nind\Remark Note that in the non-radial case, such a kind of estimates of supersolutions are extended in \cite{BiVe2} to similar estimates of the spherical minimum of $u$.
\subsection{General gradient estimates}
Here we obtain a general estimate of the gradient of a solution, non-necessarily radial, in terms of the function itself. It is based upon a combination of Bernstein and Keller-Osserman methods.

\bth{sup-2} Let $N\geq 1$ and  $q>1$. If $u\in C(\overline {B_{r_0}}\setminus\{0\})$ is a solution of $(\ref{I-1})$ in $B_{r_0}\setminus\{0\}$. 
Then for any $\gr\in (0,\frac{r_0}{2}]$ and any $x\in B_{\gr}\setminus\{0\}$
\bel{III-7}
|\nabla u(x)|\leq c_1|x|^{-\frac{1}{q-1}}+c_2\max_{B_\gr(x)}e^{\frac{u}{q}}+c_3\max_{B_\gr(x)}e^{\frac{u}{2(q-1)}},
\ee
where $c_j=c_j(N,q,M)>0$, for $j=1,2,3$.
\es
\Proof The technique is standard and we recall it for the sake of completeness. We set $z=|\nabla u|^2$, then by Schwarz inequality,
$$-\frac{1}{2}\Gd z+\frac{(\Gd u)^2}{N}+\langle\nabla \Gd u,\nabla u\rangle\leq 0.
$$ 
Then, for $\ge>0$ small enough, 
$$-\frac{1}{2}\Gd z+\frac{(Mz^{\frac q2}-e^u)^2}{N}\leq e^uz+\frac {Mq}{2}z^{\frac q2-2}|\langle \nabla z,\nabla u\rangle|\leq e^uz+\ge z^q+C_{\ge,q,M}\frac{|\nabla z|^2}{z}
$$
Using again H\"older and Young inequalities, we obtain
$$\BA{lll}\dsps
-\frac{1}{2}\Gd z+\frac{M^2z^{q}}{N}\leq 2Me^uz^{\frac q2}+e^uz+\ge z^q+C_{\ge,q,M}\frac{|\nabla z|^2}{z}\\[3mm]
\phantom{-\frac{1}{2}\Gd z+\frac{M^2z^{q}}{N}}\dsps\leq
\ge_1z^q+C_{\ge_1}e^{2u}+\ge_2z^q+C_{\ge_2}e^{q'u}+\ge z^q+C_{\ge,q,M}\frac{|\nabla z|^2}{z}
\EA$$
With the choice $\ge$, $\ge_1$ and $\ge_2$ small enough and the last inequality turns out into
\bel{III-8}
-\frac{1}{2}\Gd z+\frac{M^2z^{q}}{2N}\leq C_1e^{2u}+C_2e^{q'u}+C_{q,N,M}\frac{|\nabla z|^2}{z}
\ee
Then we use a variant of the Osserman-Keller inequality proved in \cite [Lemma 3.1]{BiVe1} and obtain $(\ref{III-7})$.\qeda\medskip

The next result which holds only if $q>2$ is a universal a priori estimate of $u$ and $\nabla u$ solution of $(\ref{I-1})$. The proof is delicate and can be found in \cite{BiVe2}
\bth{sup-3} Let $N\geq 2$ and $q>2$. If $u\in C(\overline {B_{r_0}}\setminus\{0\})$ is a solution of $(\ref{I-1})$ in $B_{r_0}\setminus\{0\}$. 
Then there exists $C>0$ depending on $N,q, M, u$ such that
\bel{III-9}
e^{u(x)}\leq \frac{C}{|x|^q}\quad\text{and }\;|\nabla u(x)|\leq \frac{C}{|x|}\quad\text{for all } x\in B_{\frac{r_0}2}\setminus\{0\}.
\ee
\es

\mysection{The dynamical system approach}
A radial solution of $(\ref{I-1})$ satisfies
\bel{II-1}
\CL^{rad}_{m,p,q}u:=-u_{rr}-\frac{N-1}{r}u_r+M|u_r|^q-e^u=0
\ee
In all the sequel we consider essentially radially symmetric solutions $u$ of  $(\ref{I-1})$ that means functions satisfying $(\ref{II-1})$. Furthermore, solutions are at least $C^3$ on their maximal interval of existence.
\subsection{Monotonicity and global properties}

\blemma{mono} Let $q>1$ and $u$ be a radial solution of $(\ref{II-1})$ on $(0,r_0)$. Then $u$ has at most one local maximum and $u_r$ keeps a constant sign near $r=0$ and 
the following dichotomy holds:\\
\nind (i) either $u_r(r)<0$ on $(0,r_0)$ and $u(r)\to\infty$ when $r\to 0$,\\
\nind (ii) or $u_r(r)<0$ on $(0,r_0)$ and $u(r)\to u_0$ when $r\to 0$,\\
\nind (iii) or $u_r(r)>0$ near $0$ and $u(r)\to u_0$ when $r\to 0$,\\
\nind (iv) or $u_r(r)>0$ near $0$ and $u(r)\to -\infty$ when $r\to 0$.
\es
\Proof The proof follows easily from the monotonicity of $u$.\qeda

\bprop{glob} Let $N\geq 2$ and  $q>1$. If $u$ is a solution of $(\ref{II-1})$ defined in  the maximal interval $(0,R)$ and decreasing near $0$, then $R=\infty$ and $u_r<0$ on $(0,\infty)$. Furthermore $u(r)\to-\infty$ and $u_r(r)\to 0$ when $r\to\infty$. As a consequence, if $u$ is positive near $0$, it has a unique zero on $(0,\infty)$.
\es
\Proof By \rlemma {mono}, $u_r<0$ on $(0,R)$. Set 
\bel{IV-2*}
\CH(r)=e^u+\frac{u_r^2}{2}.
\ee
Then
$$
\CH_r(r)=-\frac{N-1}{r}u_r^2+M|u_r|^qu_r=-\frac{N-1}{r}u_r^2-M|u_r|^{q+1}.
$$
Hence $\CH$ is decreasing on $(0,R)$. As a consequence for any $\tilde r\in (0,R)$, $u^2_r$ is bounded on $[\tilde r,R)$.  If $R<\infty$, by integration $u(r)$ is also bounded on $[\tilde r,R)$ which is impossible since $R$ is maximal. Hence $R=\infty$. Since $e^{u(r)}$ is decreasing and positive, there exists $\ell\geq 0$ such that  $e^{u(r)}\to\ell$ when $r\to\infty$. Moreover, since it is positive and decreasing $\CH(r)$ admits a limit $\gl\geq 0$ when $r\to\infty$ and therefore $u_r(r)$ shares this property. This implies that 
$\CH_r$ is integrable and therefore $u_r(r)\to 0$ when $r\to\infty$. Since $e^{u(r)}\to\ell\geq 0$, and $u_r(r)\to 0$ we obtain that $-u_{rr}(r)\to -\ell$. If $\ell>0$ we would obtain that 
$u(r)\to \ln\ell$, which is not compatible. Hence $\ell=0$ and $u(r)\to-\infty$ $r\to\infty$.\qeda\medskip

\nind {\bf Remark 3.2 bis}  More generally, if $u$ is a solution defined on a maximal interval $I_u$ containing $r_0>0$ and if $u_r(r_0)\leq 0$, then $u_r(r)<0$ for $r\in I_u\cap(r_0,\infty)=(r_0,\infty)$, 
$u(r)\to-\infty$ and $u_r(r)\to 0$ when $r\to\infty$. On the contrary, if $u_r(r_0)> 0$, then either $u$ admits a unique maximum at $r_1>r_0$ and therefore $u(r)\to-\infty$ and $u_r(r)\to 0$ when $r\to\infty$, or 
$u$ is increasing on $I_u\cap(r_0,\infty)$. In such a case two situations could occur: either $I_u\cap [r_0,\infty)=[r_0,r_1)$ and $\lim_{r\to r_1}u(r)=\infty$ ($u$ is a large solution), or $I_u\cap [r_0,\infty)=[r_0,\infty)$ and $\lim_{r\to \infty}u(r)=\infty$.
\subsection{Associated differential systems}
To the equation $(\ref{II-1})$ we associate several systems autonomous or not.
\subsubsection{Non-autonomous systems of order 2}
\blemma{L-1} Let $x$, $X$ and $\Phi$ be defined by 
\bel{II-2}
x(t)=r^2e^{u(r)}\,,\;\; X(t)=r^qe^{u(r)}\,,\;\Phi(t)=-ru_r(r)\quad\text{with }\;t=\ln r.
\ee
Then $u$ is a solution of $(\ref{II-1})$ if and only if 
\bel{II-3}
\left\{\BA{lll} 
x_t=x(2-\Phi)\\
\Phi_t=(2-N)\Phi-Me^{(2-q)t}|\Phi|^q+x,
\EA
\right.
\ee
where $x>0$. This system is also equivalent to
\bel{II-4}
\left\{\BA{lll} 
X_t=X(q-\Phi)\\
\Phi_t=(2-N)\Phi-e^{(2-q)t}\left(-M|\Phi|^q+X\right).
\EA
\right.
\ee
\es
\Proof We set 
\bel{II-5}
u(r)=U(t)\quad\text{with }t=\ln r.
\ee
Then
\bel{II-6}
-U_{tt}-(N-2)U_t-e^{2t}e^U+Me^{(2-q)t}|U_t|^q=0,
\ee
and $\Phi=-U_t=-r\frac{v_r}{v}$ with $v=e^u$. The proof follows.\qeda\medskip

\nind\Remark For the Emden equation $(\ref{I-E-C})$ the system in $(x,\Phi)$ is
$$
\left\{\BA{lll} 
x_t=x(2-\Phi)\\
\Phi_t=(2-N)\Phi+x.
\EA
\right.
$$
The equilibrium are $(0,0)$ and $(2(N-2),2)$. If $N>2$, the characteristic value of the linearisation at $(0,0)$ are $\gl_1=2-N$ with eigenvector $(0,1)$ and $\gl_2=2$ with eigenvector $(N,1)$. Hence $(0,0)$ is a saddle point. The stable trajectory at $(0,0)$ is located on $x=0$ and actually it is
$(x(t),\Phi(t))\equiv (0,ce^{(2-N)t})$. It is not admissible.  The unstable trajectory at $(0,0)$ satisfies $(x(t),\Phi(t))= (Ne^{2t},e^{2t}) (c+o(1))$ ($c>0$) when $t\to-\infty$, which corresponds to a solution $u$ 
satisfying $(u(0),u_r(0)=(\ln cN,0)$. Replacing $\ln cN$ by $u_0$, we obtain all the regular with $u(0)=u_0$.
\subsubsection{First autonomous systems of order 3}
It is well-known that a non-autonomous system of order 2 can be transformed into an autonomous system of order 3. Actually several transformations are possible.
\blemma{L2} Let $u$ be a solution of $(\ref{II-1})$. Set
\bel{II-7}
x=r^2e^u\,,\; X=r^qe^u\,,\; \Phi=-ru_r\,,\; V=r|u_r|^{q-1}\quad\text{with }\;t=\ln r.\ee
Then there holds in variable $(x,\Phi,V)$
\bel{II-8}\left\{\BA{lll}
x_t=x(2-\Gf)\\
\Phi_t=(2-N)\Phi-M|\Phi|V+x\\
\dsps
V_t=V\left(N-(N-1)q-(q-1)\left(MVsign(\Phi)-\frac{x}{\Phi}\right)\right)
\EA\right.\ee
and also with $(X,\Phi,V)$
\bel{II-8*}\left\{\BA{lll}
X_t=X(q-\Gf)\\\dsps
\Phi_t=(2-N)\Phi -V|\Phi|\left(M-\frac{X}{|\Phi|^q}\right)\\
\dsps
V_t=V\left(N-(N-1)q-(q-1)\left(M-\frac{X}{|\Phi|^q}\right)sign(\Phi)V\right)
\EA\right.\ee
\es
\Proof By a direct computation
$$r^{2-q}=e^{(2-q)t}=|\Phi|^{q-1}V,
$$
$$\Phi_t=(2-N)\Phi-M|\Phi|V+x=(2-N)\Phi+|\Phi|^{1-q}V(-M|\Phi|^q+X)$$
and, whatever is the sign of $\Phi$,
$$\BA{lll}\dsps
\frac {V_t}V=(2-q)+(q-1)\frac {\Phi_t}\Phi=2-q+(q-1)\left(2-N+\frac{-M|\Phi| V+x}{\Phi}\right)\\
\phantom{\dsps\frac {V_t}V}\dsps=q-N(q-1)+(q-1)\frac{-M|\Phi| V+x}{\Phi}\\
\phantom{\dsps\frac {V_t}V}\dsps=q-N(q-1)+(q-1)\left(-MVsign (\Phi)+\frac x\Phi\right)\\
\phantom{\dsps\frac {V_t}V}\dsps=q-N(q-1)+(q-1)\left(-M+\frac X{|\Phi|^q}\right)Vsign (\Phi),
\EA$$
which leads to $(\ref{II-7})$ and $(\ref{II-8*})$.\qeda\medskip

We also introduce another system of order 3 in the variables $(x,\Phi,\Gth)$ where $\Gth(t)=e^{(2-q)t}$, 
\bel{II-9}\left\{\BA{lll}
x_t=x(2-\Gf)\\
\Phi_t=x+(2-N)\Phi-M|\Phi|^q\Gth\\
\dsps
\Gth_t=(2-q)\Gth.
\EA\right.\ee
This system will be interesting when $t\to-\infty$ (i.e. singularity in $x$) if $1<q<2$ and when $t\to\infty$ if $q>2$ since in that in these two cases $\Gth(t)\to 0$ when $t\to-\infty$ or $t\to\infty$ according $r\to 0$ or $r\to\infty$. It admits two equilibria, $(0,0,0)$ and $(2(N-2),2,0)$.
\subsubsection{Reduction to an autonomous quadratic system of order 3}
Using a suitable change of variable, the equation is transformed into a remarkable  quadratic system of Lotka-Volterra type as the next lemma shows it.

\blemma{L3} Let $u$ be a solution of $(\ref{II-1})$. At any point where $u_r(r)\neq 0$ define
\bel{II-11}
Z=-\frac{re^u}{u_r}\,,\; V=r|u_r|^{q-1}\,,\; \Phi=-ru_r
\quad\text{with }\;t=\ln r.\ee
Then there holds 
\bel{II-12}\left\{\BA{lll}
Z_t=Z(N-\Phi+sMV-Z)\\
\dsps
V_t=V\left(N-(N-1)q+(q-1)(Z-sMV)\right)\\
\Phi_t=\Phi\left(2-N+Z-sMV\right)
\EA\right.\ee
where $s=sign(\Phi)$.
\es
\Proof We start from $(\ref{II-8})$ and define $Z=\frac{x}{V}=-\frac{re^u}{u_r}$. Then
$$\frac {Z_t}{Z}=\frac{x_t}{x}-\frac{\Phi_t}{\Phi}=2-\Phi +N-2+M\frac{|\Phi|}{\Phi}-Z,
$$
which implies $(\ref{II-12})$.\qeda\medskip

\nind\Remark For $q\neq 2$ the system $(\ref{II-11})$ admits four equilibria
$$O=(0,0,0)\,,\; Q_0=(N-2,0,2)\,,\; N_0=(N,0,0)\,,\; P_0=(0,s\tfrac{(N-1)q-N)}{q-1},0).
$$
The linearised system at $Q_0$ is
\bel{II-13}\left\{\BA{lll}
Z_t=(N-2)(N-2-Z+MV-\Phi)\\
\dsps
V_t=\left(2-q\right)V\\
\Phi_t=\left(2-N+Z-MV\right).
\EA\right.\ee
The characteristic polynomial is
$$P(\gl)=(2-q-\gl)(\gl^2+(N-2)\gl+2(N-2)).
$$
One root is $2-q$, the two other roots have negative real part. Hence if $1<q<2$ there exist a 2-dimensional manifold of trajectories converging to $Q_0$ when $t\to\infty$ and an unstable trajectory issued from $Q_0$. This unstable trajectory corresponds to a solution $u$ of $(\ref{I-1})$ satisfying $\dsps\lim_{r\to 0}r^2e^{u(r)}=2(N-2)$. If $q>2$, $Q_0$ is a sink which attracts all local trajectories. \smallskip

\nind The linearised system at $N_0$ is 
\bel{II-14}\left\{\BA{lll}
Z_t=N\left(-\Phi-MV-Z+N\right)\\
\dsps
V_t=qV\\
\Phi_t=2\Phi.
\EA\right.\ee
The characteristic polynomial is
$$P(\gl)=-(q-\gl)(2-\gl)(N+\gl).
$$
It can be checked that the trajectories converging to $N_0$ at $-\infty$ correspond to the regular trajectories.\smallskip

\nind Finally the trajectories converging to $P_0$ at $-\infty$ with $q>2$ correspond to the solutions of Hamilton-Jacobi type as it is shown in the proof of \rth{q>2-beta*}.

\mysection{Isolated singularities when $1<q<2$}

\subsection{Singular behaviour}

In this section we use a perturbation argument due to \cite[Proposition 4.1]{LoRy} that we recall below.

\bprop{LR} Let $h:\BBR_+\ti \BBR^N\to\BBR^N$ be a Carath\'eodory function. Assume that there exists locally continuous  function $h^*:\BBR^N\to\BBR^N$ such that for all 
compact set $C\subset\BBR^N$ and all $\ge>0$ there exists $T=T(\ge,C)>0$ such that 
$${\rm {ess}}\sup_{\!\!\!\!\!\!x\in C}\sup_{t\geq T}|h(\gt,x)-h^*(x)|\leq \ge.
$$
If $x(t)$ is a bounded solution of $x_\gt=h(\gt,x)$ on $\BBR_+$ such that $x(0)=x_0$,
then the omega-limit set of the positive trajectory of $x$ is a no-empty connected compact set of $\BBR^N$ which is invariant under the flow of the equation $x_\gt=h^*(x)$.
\es

\nind{\it Mutatis mutandis} a similar result holds if $\BBR_+$ is replaced by $\BBR_-$ and omega-limit set by alpha-limit set.\\

 Our first result is a complete description of the behaviour of any solution near $0$.

\bth{compsub} Let $1<q<2$ and $N\geq 3$. If $u$ is any radial solution of $(\ref{I-1})$ in $B_{r_0}\setminus\{0\}$ there holds\smallskip

\nind (i) either 
\bel{IV-3*}\lim_{r\to 0}r^2e^{u(r)}=2(N-2)\,,\;\text{ that is } \;\lim_{r\to 0}\left(u(r)-2\ln\frac 1r\right)=\ln 2(N-2).
\ee 
\nind (ii) or $u$ is a regular solution, i.e. there exists $u_0\in\BBR$ such that $\dsps\lim_{r\to 0}u(r)=u_0$ and $\dsps\lim_{r\to 0}u_r(r)=0$,\smallskip

\nind (iii) or $\dsps\lim_{r\to 0}u(r)=-\infty$ and
\bel{IV-4*}
\lim_{r\to 0}r^{N-2}u(r)=\gg<0\quad\text{when }\, 1<q<\frac{N}{N-1},
\ee
\bel{IV-5log}
\lim_{r\to 0}r^{N-2}|\ln r|^{N-1}u(r)=-\frac{1}{N-2}\left(\frac{N-1}{M}\right)^{N-1}\quad\text{when }\, q=\frac{N}{N-1},
\ee
\bel{IV-6}
\lim_{r\to 0}r^{\frac{2-q}{q-1}}u(r)=-\xi_M:=-\frac{q-1}{2-q}\left(\frac{(N-1)q-N}{M(q-1)}\right)^{\frac{1}{q-1}}\quad\text{when }\, \frac{N}{N-1}<q<2.
\ee

\es 
\Proof By \rth{sup-1} we have that $r^2e^{u(r)}\leq C=C(N,p,q,M)>0$ when $0<r\leq \frac {r_0}2$. We consider the system $(\ref{II-3})$ in $(x,\Phi)=(r^2e^{u(.)},-ru_r)$ in the variable $t=\ln r$, that is 
\bel{IV-7}\left\{\BA{lll}
x_t=x(2-\Phi)\\
\Phi_t=x+(2-N)\Phi
-Me^{(2-q)t}|\Phi|^q.\EA\right.\ee
{\it Step 1: }We assume that $u$ is decreasing near $r=0$. Then $\Phi>0$ and $x(r)\leq C$ for $0<r\leq \frac{r_0}2$.\\
 If $\Phi$ is unbounded on $(-\infty, \tilde T]$ with $\tilde T=\ln\tilde r$, we encounter two possibilities:\\
- either $\Phi$ is monotone as $t\to -\infty$ and $\Phi(t)\to \infty$. Then for any $A>0$ we have that $x_t(t)\leq -Ax(t)$ for $t\leq t_A\leq \tilde T$, which implies that $t\mapsto e^{At}x(t)$ is decreasing on 
$(-\infty,t_A]$, thus $x(t)\geq e^{A(t_A-t)}x(t_A)$ on this interval, an inequality which contradicts the boundedness of $x(t)$. \\
- or $\Phi$ is not monotone and thus there exists a sequence $\{t_n\}$ tending to $-\infty$ such as $\Phi(t_n)$ is a local maximum of $\Phi$, and the sequence $\{\Phi(t_n)\}$ tends to $\infty$ when $n\to\infty$. From $(\ref{IV-7})$ we have that $x(t_n)=(N-2)\Phi(t_n)
+Me^{(2-q)t_n}|\Phi(t_n)|^q\to\infty$.  This contradicts the boundedness of $x$. Hence $\Phi$ is bounded. therefore there exists $\gk>0$ such that 
$\sup\{|x(t)|,|\Phi(t)|\}\leq \gk$  for $t\leq \ln\tilde T$.  Since $q<2$ this system is an exponential perturbation as $t\to-\infty$ of the sytem  
\bel{IV-8*}\left\{\BA{lll}
x_t=x(2-\Phi)\\
\Phi_t=x+(2-N)\Phi,\EA\right.\ee
in the sense of \rprop{LR}, which is the system associated to the Chandrasekhar-Emden equation $-\Gd u=e^u$. This system admits the equilibria $(2(N-2),2)$ and $(0,0)$ which both are hyperbolic. The point $(2(N-2),2)$ is a sink while $(0,0)$ is a saddle point. Therefore $(2(N-2),2)$ is repulsive when $t\to-\infty$ and the only solution in its neighbourhood is the constant solution $(2(N-2),2)$ which corresponds to the explicit solution of the Chandrasekhar-Emden equation, $u(r)=\ln\frac{2(N-2)}{r^2}$. To the stable trajectory of $(\ref{IV-8*})$ converging to $(0,0)$ corresponds the regular solution of the same equation.\\
By  \rprop{LR} any bounded solution of $(\ref{IV-7})$ admits a limit set at $-\infty$ which is invariant under the flow of $(\ref{IV-8*})$, actually, the only possibilities are $(2(N-2),2)$ and $(0,0)$. Then either $(x(t),\Phi(t))\to (2(N-2),2)$ or $(x(t),\Phi(t))\to (0,0)$ when $t\to-\infty$.\\
 In the first case $u(r)$ satisfies $(\ref{IV-3*})$. In the second case we introduce the system $(\ref{II-9})$ with 
$\Gth(t)=e^{(2-q)t}$ and we are in the situation where $(x(t,\Phi(t),\Gth(t))\to (0,0,0)$ when $t\to-\infty$. We recall the system  $(\ref{II-9})$ in $(x,\Phi,\Gth)$ 
with $\Gth(t)=e^{(2-q)t}$, 
\bel{IV-7-2}\left\{\BA{lll}
x_t=x(2-\Gf)\\
\Phi_t=x+(2-N)\Phi-M|\Phi|^q\Gth\\
\dsps
\Gth_t=(2-q)\Gth.
\EA\right.\ee
The associated linearised system  at $(0,0,0)$ is 
\bel{IV-9}\left\{\BA{lll}
\;x_t=2x\\
\,\Phi_t=x+(2-N)\Phi\\
\,\Gth_t=(2-q)\Gth.
\EA\right.\ee
The eigenvalues are $\gl_1=2>0$, $\gl_2=2-N<0$ and $\gl_3=2-q>0$. Hence there exists  a 2-dimensional unstable manifold $\CM_s$ of trajectories issued of $(0,0,0)$ as $t\to-\infty$; $\CM_s$ is relative to the eigenvalues 
$\gl_1$ and $\gl_3$ with corresponding eigenvectors $\gw_1=(N,1,0)$ and $\gw_3=(0,0,1)$ which generate the tangent $2$-plane at this point. In the manifold $\CM$ there exists one trajectory $\CT_f$, called the fast trajectory,  
associated to $\gl_1$ and admitting the vector $\gw_1$ for tangent vector at this point. However this trajectory is located in the plane $\Gth=0$ since the point $(0,0)$ of the restriction of $(\ref{IV-7-2})$ to this plane is a saddle point, hence this trajectory is also the unstable trajectory of $(0,0)$. Therefore this trajectory is not admissible. 
Hence our trajectory is associated to $\gl_3$ with tangent  vector  $\gw_3$ at $(0,0,0)$.
Along this trajectory there holds $x(t)=o(e^{(2-q)t})=o(\Gth(t))$ and $\Phi(t)=o(e^{(2-q)t})$. Then 
$r|u_r(r)|=o(r^{2-q})$ when $r\to 0$, and thus $u(r)$ has a finite limit $u_0$ at $r=0$, which implies $x(t)=e^{u_0}e^{2t}(1+o(1))$ as $t\to-\infty$. It follows from the equation $(\ref{I-1})$ and elliptic equation regularity that 
$u$ can be extended as a $C^2$ solution, thus $u_r(0)=0$ and $u$ is a regular solution.\smallskip

\nind{\it Step 2: }We assume that $u$ is increasing near $r=0$. Then either $u(r)$ has a finite limit $u_0$ when $r\to 0$ or $u(r)\to-\infty$. Moreover $u_r$ is monotone near $0$: indeed  at any point $\tilde r$ where $u_{rr}(\tilde r)=0$ the following identity
$$u_{rrr}(\tilde r)=\left(\frac{N-1}{\tilde r}-e^{u(\tilde r)}\right)u_{r}(\tilde r)>0.
$$
Then \\
- either $u_r(r)\to 0$ but in that case $u(r)\to u_0$ at $0$ and $u$ is concave near $0$, hence it is decreasing, contradiction,\\
- or $u_r$ has a positive limit $c_0$ at $0$, then again  $u(r)\to u_0$ but from the equation $u_{rr}$ is not integrable at $0$, contradiction,\\ 
- or $u_r$ tends to $\infty$. In that case $e^{u(r)}=o(u^q_r(r))$. In that case equation $(\ref{II-1})$ can be written under the form 
\bel{IV-10}-u_{rr}-\frac{N-1}{r}u_r+\tilde M(r)u^q_r=0
\ee
where $\tilde M(r)=M-(u_r(r)^{-q})e^{u(r)}\to M$ when $r\to 0$. Equation $(\ref{IV-10})$ is explicitely integrable and we
have
\bel{IV-11}\BA{lll}\dsps r^{-(N-1)(q-1)}X(r)=r_0^{-(N-1)(q-1)}X(r_0)+(q-1)\int_r^{r_0}s^{-(N-1)(q-1)}\tilde M(s)ds.
\EA\ee
where $X(r)=u_r^{1-q}$. By performing a direct integration of $(\ref{IV-11})$ we obtain (iii).\qeda
\subsection{Existence of singular solutions of Emden-Chandrasekhar type}

\bth{posi} Let $1<q<2$ and $N\geq 3$. Then there exists a unique radial solution $u_\gw$ of $(\ref{II-1})$ defined on $(0,\infty)$ such that 
\bel{Ex-1}
\lim_{r\to 0}r^2e^{u(r)}=2(N-2)\quad\text{and }\,\lim_{r\to 0}ru_r(r)=-2.
\ee
\es
\Proof We still consider the systems $(\ref{IV-7})$ and $(\ref{IV-7-2})$ relative to $(x,\Phi,\Gth)$. Here we prove the existence of a unique trajectory $\CT$ of $(\ref{IV-7-2})$ converging to the point $P_0=(2(N-2),2,0)$ as $t\to-\infty$. The linearised system at $P_0$ with $x=2(N-2)+\bar x$ and $\Phi=2+\bar\Phi$ is 
 \bel{Ex-4}\left\{\BA{lll}
\bar x_t=-2(N-2)\bar \Phi\\
\bar \Phi_t=\bar x+(2-N)\bar \Phi-M2^q\bar \Gth\\
\Gth_t=(2-q)\Gth.
\EA\right.
\ee
Set 
$$A=\left(\BA{llll} 0\;&-2(N-2) \!\!\!\!\!\!&&0\\
1\; &\;\;\;2-N \!\!\!\!\!\!\!\!&&-2^qM\\
0\; &\;\;\;\;\;\;0 \!\!\!\!\!\!&&2-q \EA\!\!\right)
$$
The characteristic polynomial associated is
 \bel{Ex-VP}\det(A-\gl I):=P(\gl)=(2-q-\gl)\left(\gl^2+(N-2)\gl+2(N-2)\right).
\ee

The corresponding eigenvalues are $\gl_1=2-q>0$ and $\gl_2, \gl_3$ which are real and negative if $N\geq 10$ and complex with negative real part if $3\leq N\leq 9$.  Hence we have the standard decomposition
$$\BBR^3=\ker (A-\gl_1I)\oplus H
$$
where $H$ is either $\ker (A-\gl_2I)\oplus \ker (A-\gl_3I)$ if $\gl_2\neq \gl_3$ are real, or  $\ker (A-\gl_2I)^2$  if $\gl_2=\gl_3$ (necessarily real), or 
$\CR e\left(\ker (A-\gl_2I)\oplus \ker (A-\gl_3I)\right)$ if $\gl_2\neq \gl_3$ are not real but conjugate. Then $\ker (A-\gl_1I)=$span$ \,\{\gw_1\}$ where $\gw_1$ has for components
$$\gw_1=(2(N-2),q-2,2^qMf(q))\,\text{ where }\;f(q)=q^2-(N+2)q +4(N-1).
$$
We prove that that $f(q)\neq 0$. This is clear if $N<10$. If $N\geq 10$ it admits two positive roots $q_1<q_2$, with $2<q_2$. Since $f(2)>0$, then $2<q_1$. Therefore 
$f(q)>0$. Then
 there exists a unique trajectory $\CT^*$ associated to $(x^*,\Phi^*,\Gth^*)$ such that 
$\dsps \lim_{t\to-\infty}(x^*(t),\Phi^*(t),\Gth^*(t))=(2(N-2),2,0)$ with tangent vector at the trajectory at $P_0$ is colinear to $\gw_1$ and such that $\Gth^*(t)>0$ as $t\to\infty$.  To this trajectory is associated a solution $u_\gw$ of $(\ref{II-1})$ satisfies 
$(\ref{Ex-1})$. It is decreasing; by \rprop{glob} it is defined on $(0,\infty)$ and it satisfies $\dsps \lim_{r\to\infty}u(r)=-\infty$ and $\dsps \lim_{r\to\infty}u_r(r)=0$. Any other solution corresponding to the same trajectory is just a time shift of $(x^*(t),\Phi^*(t),\Gth^*(t)$ and it corresponds to the same function $u_\gw$. \qeda\medskip

\subsection{Existence of singular solutions of Hamilton-Jacobi type}
Next we show the existence of solutions satisfying $(\ref{IV-6})$.
\bth{nega} Let $\frac N{N-1}<q<2$. Then there  exists infinitely many radial solutions $u$ of $(\ref{II-1})$ defined on $(0,\infty)$ such that 
 \bel{Ex-5}\BA{lll}\dsps
\lim_{r\to 0}r^\gb u(r)=-\xi_M,
\EA\ee
where $\gb=\frac {2-q}{q-1}$.
\es
\Proof We set $U=-u$, then 
$$-\Gd U-M|\nabla U|^q+e^{-U}=0,
$$
and we put
 \bel{Ex-6}
U(r)=r^{-\gb}\xi(t)\,,\;U_r(r)=-r^{-\frac1{q-1}}\eta(t)=-r^{-\gb-1}\eta(t)\;\text{ with }\,t=\ln r.
\ee
We are led to the system
 \bel{Ex-7}\left\{\BA{lll}
\xi_t=\beta\xi-\eta\\
\eta_t=-\gk\eta+M\eta^q-e^{\frac{qt}{q-1}}e^{-e^{-\gb t}\xi(t)},
\EA\right.\ee
with $\gk=\frac{(N-1)q-N}{q-1}>0$. The couple $(\xi,\eta)$ is a solution if $(\xi(t),\eta(t))\to(\xi_M,\gb\xi_M)$ when $t\to-\infty$. Set 
$\bar \xi=\xi-\xi_M$ and $\bar \eta=\eta-\gb\xi_M$. Then 
 \bel{Ex-8}\left\{\BA{lll}
\bar\xi_t=\beta\bar\xi-\bar\eta\\
\bar\eta_t=\gk(q-1)\bar\eta+F(\bar\eta)-e^{\frac{qt}{q-1}}e^{-e^{-\gb t}(\bar\xi +\xi_M)},
\EA\right.\ee
where $F(\bar\eta)=O(|\bar\eta|^{2})$. Let $0<\gth<4\gk (2-q)$, then there exists $\gm:=\mu(\gth,\gk)>0$ such that  
$$\theta\beta x^2-\theta xy+\gk(q-1) y^2\geq \gm(\theta x^2+y^2).
$$
Hence 
 \bel{Ex-9}\frac 12\frac{d}{dt}(\theta \bar\xi^2+\bar\eta^2)\geq\gm (\theta \bar\xi^2+\bar\eta^2)+\bar\eta\left(F(\bar\eta)-e^{\frac{qt}{q-1}}e^{-e^{-\gb t}(\bar\xi +\xi_M)}\right)
\ee
There exists $(x_0,y_0)$  such that 
$$\BA{lll}
(i)&\dsps|x_0|\leq\frac 14\xi_M\\[2mm]
\dsps
(ii)&\dsps|y_0|\leq \min \{\frac \gb4\xi_M,c\}\Longrightarrow |F(y_0)|\leq \frac\gm 2|y_0|^3\leq \frac\gm 2|y_0|^2
\EA$$
and if  we choose $(\bar\xi_0,\bar\eta_0)$  such that  $|\bar\xi_0|\leq |x_0|$ and $|\bar\eta_0|\leq |y_0|$. We denote by $(\bar\xi(t),\bar\eta(t))_{t\leq t_0}$ the solution of $(\ref{Ex-8})$ with initial data 
$(\bar\xi(t_0),\bar\eta(t_0))=(\bar\xi_0,\bar\eta_0)$. As long as $|\bar\xi(t)|\leq |x_0|$ and $|\bar\eta(t)|\leq |y_0|$ we have
$$|\eta|e^{\frac{qt}{q-1}}e^{-e^{-\gb t}(\bar\xi +\xi_M)}\leq \frac\gm4\eta^2+\frac 4\gm e^{\frac{2qt}{q-1}}e^{-2e^{-\gb t}\frac{\xi_M}{2}}
$$
and 
$$\frac 12\frac{d}{dt}(\theta \bar\xi^2+\bar\eta^2)\geq\frac\gm2(\theta \bar\xi^2+\bar\eta^2)-\frac {c(t)}2,
$$
where $C(t)$ is a positive function which satisfies 
$$\lim_{t\to-\infty}e^{-at}c(t)=0\quad\text{for all }a>0.
$$
This implies
 \bel{Ex-10}
 \theta \bar\xi(t)^2+\bar\eta(t)^2\leq e^{\gm(t-t_0)}  (\theta \bar\xi_0^2+\bar\eta_0^2)
 +e^{\gm t}\int_t^{t_0}e^{-\gm s }c(s) ds\quad\text{for }t\leq t_0.
\ee
As long as $|\bar\xi(t)|\leq |x_0|$ and $|\bar\eta(t)|\leq |y_0|$, the above inequality holds. If we take $\theta \bar\xi_0^2+\eta_0^2\leq \min\{\theta x_0^2,y_0^2\}$, inequality 
$(\ref{Ex-10})$ holds for all $t\leq t_0$. Hence $(\bar \xi(t),\bar\eta(t))\to (0,0)$ when $t\to-\infty$, which implies 
 \bel{Ex-11}
\lim_{t\to-\infty}(\xi(t),\eta(t))=(\xi_M,\gb\xi_M).
\ee
If we set $r_0=e^{t_0}$ and $u(r)=-r^{-\gb}\left(\xi_M+\bar\xi(\ln r)\right)$, then $u$ is a solution of $(\ref{II-1})$ in $(0,r_0]$ which satisfies $(\ref{Ex-5})$. Such a solution can be extended to 
$(0,\infty)$ by \rprop{glob} and the choice of its data at $r=r_0$ has for unique restriction $u(r_0)$ and $u_r(r_0)$ corresponding to $(\bar\xi_0,\bar\eta_0)$.\qeda\medskip

Next we show the existence of solutions of $(\ref{II-1})$ satisfying $(\ref{IV-4*})$ by a fixed point method.

\bth{dirac} Let $N\geq 3$ and $1<q<\frac N{N-1}$. Then there exists $\gr_0>0$ and $k_0>0$ such that for $0<\gr\leq \gr_0$ and $-k_0<\gg\leq 0$ there exists a radial function $u_\gg$ satisfying
 \bel{Ex-12}\BA{lll}
-\Gd u_\gg+M|\nabla u_\gg|^q-e^{u_\gg}=c_N\gg\gd_0\qquad&\text{in }\CD'(B_\gr)\\
\phantom{-\Gd +M|\nabla u_\gg|^q-e^{u_\gg}}
u_\gg=0&\text{on }\prt B_\gr.
\EA\ee
Furthermore
 \bel{Ex-13*}
\lim_{r\to 0}r^{N-2}u_\gg(r)=\gg.
\ee
\es
\Proof We look for a radial solution $u$satisfying
 \bel{Ex-14}
\lim_{r\to 0}r^{N-2}u(r)=\gg \,\text{and }\;\lim_{r\to 0}r^{N-1}u(r)=(1-N)\gg.
\ee
If such a solution exists $e^u$ and $|\nabla u|^q$ are integrable and $(\ref{Ex-12})$ holds. 
The function $U=-u$ has to satisfy (note that $-\gg>0$)
 \bel{Ex-13}\BA{lll}\dsps
-U_{rr}-\frac{N-1}{r}U_r-M|\nabla U|^q+e^{-U}=0\qquad\text{in }(0,\gr)\\[2mm]
\phantom{,,,,,,,,,,,,,,,,,,,,,}\dsps\!\lim_{r\to 0}r^{N-2}U(r)=-\gg\\[2mm]
\phantom{,,,,,,,,,,,,,,,,,,,,,}\dsps\!\!\lim_{r\to 0}r^{N-1}U_r(r)=(N-1)\gg\\[2mm]
\phantom{,,,,,,,,,,,,,,,,,,,,,,,,,,,,,,,}\!\!U(\gr)=0.
\EA\ee
Hence the function $V(r)=U_r(r)$ satisfies
 \bel{Ex-15}
V(r)=(N-1)\gg r^{1-N}-r^{1-N}\int_0^r\left(M|V|^q-e^{-U}\right)s^{N-1}ds.
\ee
We combine this with 
 \bel{Ex-16}
U(r)=-\int_r^\gr V(s)ds,
\ee
and define the operator $(U,V)\mapsto K(U,V)=(K_1(U,V),K_2(U,V))$ with
 \bel{Ex-17}\BA{lll}\dsps
K_1(U,V)(r)=-\int_r^\gr V(s)ds\\[2mm]
\dsps
K_2(U,V)(r)=(N-1)\gg r^{1-N}-r^{1-N}\int_0^r\left(M|V|^q-e^{-|U|}\right)s^{N-1}ds.
\EA\ee
We define $K$ on the subspace $\CK$ of $C((0,\gr])\times C((0,\gr])$ of functions $W=(U,V)$ which satisfy 
$$\norm{W}_\CK=\norm{(U,V)}_\CK=\max\left\{\gs\sup_{0<r<\gr}r^{N-2}|U(r)|,\sup_{0<r<\gr}r^{N-1}|V(r)|\right\}:=\max\left\{\gs N_1(U),N_2(V)\right\}<\infty,
$$
where $0<\gs<1$.\smallskip

\nind{\it Step 1: Lipschitz estimate.} We have
$$\BA{lll}\dsps \norm{K(U_1,V_1)-K(U_2,V_2)}_\CK=\max\left\{\gs\sup_{0<r<\gr}r^{N-2}\left|\int_r^\gr(V_1-V_2)ds\right|,\right.\\[4mm]
\phantom{----------------}\dsps\left.\sup_{0<r<\gr}\left|\int_0^r\left(M\left(|V_1|^q-|V_2|^q\right)-\left(e^{-|U_1|}-e^{-|U_2|}\right)\right)s^{N-1}ds\right|\right\}\\[4mm]
\phantom{\dsps \norm{K(U_1,V_1)-K(U_2,V_2)}_\CK}=\max\left\{I_1,I_2\right\}.
\EA$$
Since for $r<s<\gr$
$$|(V_1-V_2)(s)|\leq s^{1-N}\sup_{0<r<\gr}r^{N-1}|(V_1-V_2)(r)|,
$$
we have that
$$I_1\leq \frac\gs{N-2}\sup_{0<r<\gr}r^{N-1}|(V_1-V_2)(r)|=\frac{\gs}{N-2}N_2(V_1-V_2).
$$
Concerning $I_2$, we have
$$\left|(e^{-|U_1|}-e^{-|U_2|}\right|\leq ||U_1|-|U_2||\leq |U_1-U_2|,
$$
hence
$$\BA{lll}\dsps\sup_{0<r<\gr}\left|\int_0^r\left(e^{-|U_1|}-e^{-|U_2|}\right)s^{N-1}ds\right|\leq \left(\sup_{0<r<\gr}r^{N-2}|U_1(r)-U_2(r)|\right)\int_0^\gr sds\\[4mm]
\phantom{\dsps\sup_{0<r<\gr}\left|\int_0^r\left(e^{-|U_1|}-e^{-|U_2|}\right)s^{N-1}ds\right|}\dsps\leq \frac {\gr^2}2\sup_{0<r<\gr}r^{N-2}|U_1(r)-U_2(r)|\\[4mm]
\phantom{\dsps\sup_{0<r<\gr}\left|\int_0^r\left(e^{-|U_1|}-e^{-|U_2|}\right)s^{N-1}ds\right|}\dsps\leq \frac {\gr^2}2N_1(U_1-U_2),
\EA$$
and
$$\BA{lll}\dsps \left||V_1|^q-|V_2|^q\right|(s)\leq q\sup\{|V_1(s)|^{q-1},|V_2(s)|^{q-1}\}|V_1-V_2|(s)\\[4mm]
\phantom{\dsps \left||V_1|^q-|V_2|^q\right|(s)}\dsps \leq q\sup_{0<r<\gr}\max\left\{|r^{N-1}V_1(r)|^{q-1},|r^{N-1}V_2(r)|^{q-1}\right\}s^{(q-1)(1-N)}|V_1-V_2|(s),
\EA$$
$$\BA{lll}\dsps \sup_{0<r<\gr}\left|\int_0^rm\left(|V_1|^q-|V_2|^q\right)s^{N-1}ds\right|\\[4mm]
\phantom{-----}
\dsps\leq \sup_{0<r<\gr}\frac{mq\max\left\{|r^{N-1}V_1(r)|^{q-1},|r^{N-1}V_2(r)|^{q-1}\right\}}{N-q(N-1)}\sup_{0<r<\gr}r^{N-1}|V_1(r)-V_2(r)|\\[4mm]
\phantom{-----}
\dsps\leq \sup_{0<r<\gr}\frac{mq\max\left\{|r^{N-1}V_1(r)|^{q-1},|r^{N-1}V_2(r)|^{q-1}\right\}}{N-q(N-1)}N_2(V_1-V_2).
\EA$$
Finally,
 \bel{Ex-18}\BA{lll}\dsps 
\norm{K(U_1,V_1)-K(U_2,V_2)}_\CK\leq \max\left\{ \frac \gs{N-2}N_2(V_1-V_2),\right.\\[3mm]
\phantom{\norm{K(U_1,V_1)---,,}}\dsps\left.\frac {\gr^2}2N_1(U_1-U_2)+\frac{Mq\max\{N^{q-1}_2(V_1),N^{q-1}_2(V_2)\}}{N-q(N-1)}N_2(V_1-V_2)\right\}.
\EA\ee
{\it Step 2: Bounds on the mapping $K$.} We still have to estimate the terms $N^{q-1}_2(V_j)$ and for such a task we have to find a ball $B^{\CK}$ in $C((0,\gr])\times C((0,\gr])$ endowed with the norm 
$||\,.\,||_\CK$ which is invariant under $K$. If $(U,V)\in B^{\CK}=B^{\CK}_R((0,0))$ we have by assumption
\bel{Ex-19}|U(r)|\leq\frac R\gs r^{2-N} \,\text{ and }\;|V(r)|\leq Rr^{1-N}\,\text{ for all }\;0<r\leq\gr.
\ee
Then
\bel{Ex-20}N_1(K_1(U,V)\leq \sup_{0<r\leq \gr}r^{N-2}\int_r^\gr Rs^{1-N}ds\leq  \frac R {N-2}.
\ee
Since by $(\ref{Ex-19})$
$$\int_0^r\left|M|V|^q-e^{-|U|}\right|s^{N-1}ds\leq \frac{r^N}{N}+\frac{MR^qr^{N-q(N-1)}}{N-q(N-1)}
$$
we obtain
 \bel{Ex-21}\BA{lll}\dsps 
N_2(K_2(U,V))\leq \frac {\gr^2}2N_1(U)+(N-2)|\gg|+\frac{\gr^N}{N}+\frac{MR^q\gr^{N-q(N-1)}}{N-q(N-1)},
\EA\ee
and finally
 \bel{Ex-22}
\norm{K(U,V)}_\CK\leq \max\left\{\frac{\gs R} {N-2},\frac{\gr^2R}{2\gs}+(N-2)|\gg|+\frac{\gr^N}{N}+\frac{MR^q\gr^{N-q(N-1)}}{N-q(N-1)}\right\}
\ee
We fix $\gs=\frac 34$. Therefore for any $R>0$ there exist $0<\gr_0<1$ and $k_0>0$ such that for $0\leq|\gg|\leq k_0$ and $0<\gr\leq\gr_0$ there holds
 \bel{Ex-23*}
\norm{(U,V)}_\CK\leq R\Longrightarrow \norm{K(U,V)}_\CK\leq R.\ee
{\it Step 3: End of the proof.} With the choice of $\gs$, we have $\frac \gs{N-2}=\frac{3}{4(N-2)}\leq \frac 34$. Furthermore, if 
$\norm{(U,V)}_\CK\leq R$ we have also $\dsps\sup_{0<r\leq \gr}r^{N-1}|V(r)|\leq R$, hence $N^{q-1}_2(V)\leq R^{q-1}$. Combining this fact with estimate $(\ref{Ex-18})$ 
we obtain that if $(U_1,V_1)$ and $(U_2,V_2)$ are in $B_R^{\CK}$ there holds
if 
 \bel{Ex-23}\BA{lll}\dsps 
\norm{K(U_1,V_1)-K(U_2,V_2)}_\CK\leq \max\left\{ \frac 3{4}N_2(V_1-V_2),\right.\\[3mm]
\phantom{\norm{K(U_1,V_1)-------,}}\dsps\left.\frac {\gr^2}{2}N_1(U_1-U_2)+\frac{MqR^{q-1}}{N-q(N-1)}N_2(V_1-V_2)\right\}
\EA\ee
Up to reducing the value of $R$, always with $R\leq 1$, we can assume that $\frac{mqR^{q-1}}{N-q(N-1)}\leq\frac 14$. 
A straightforward verification shows that 
 \bel{Ex-24}\BA{lll}\dsps \max\left\{ \frac 34N_2(V_1-V_2),\frac {\gr^2}{2}N_1(U_1-U_2)+\frac{1}{4}N_2(V_1-V_2)\right\}\leq \max\left\{\frac34,\gr^2\right\}\norm{U-V}_\CK.
\EA\ee
Since $\gr<1$, the mapping $K$ admits a fixed point $(U,V)$. Hence $U$ satisfies
$$\BA{lll}\dsps
-U_{rr}-\frac{N-1}{r}U_r-M|\nabla U|^q+e^{-|U|}=0\quad\text{in }(0,\gr)\\
\dsps \phantom{----------}
\lim_{r\to 0}r^{N-2}U(r)=-\gg
 \EA$$
 and it vanishes on $r=\gr$. Since $-\gg>0$ and $\lim_{r\to 0}r^{N-2}U(r)=-\gg$ there exists  $\gr_1\in (0,\gr]$ such that $U(r)>0$ for $0<r\leq \gr_1$. Hence $e^{-|U|}=e^{-U}$. Thus 
 $u=-U$ satisfies $(\ref{Ex-12})$ in $B_{\gr_1}$. By \rprop {glob} $u$ can be extended as a solution of 
  \bel{Ex-25}
 -\Gd u+M|\nabla u|^q-e^{u}=c_N\gg\gd_0\qquad\text{in }\CD'(\BBR^N).
 \ee
 Note that the existence of solutions to $(\ref{Ex-25})$ is natural as soon as it is proved that $|\nabla u|^q$ is integrable by using the classical result of \cite{BrLi}.
 \qeda
\mysection{Isolated singularities when $q>2$}
We first denote by $U_{eik}$ the solution of the eikonal equation in $\BBR^N\setminus\{0\}$
\bel{IV-3}
M|u_r|^q-e^u=0.
\ee
Its expression is
\bel{IV-4'}
U_{eik}(r)=\ln\frac{Mq^q}{r^q}=-q\ln r+\Gl_{N,q,M}\quad\text{with }\Gl_{N,q,M}=\ln Mq^q .
\ee
Notice that if $N=2$, $U_{eik}$ is a solution of $(\ref{II-1})$ in $\BBR^2\setminus\{0\}$. If $N\geq 2$ it is a supersolution. \\
Our first key result is the description of singularities at $r=0$ when $q>2$.
\subsection{Singular behaviour}

\bth{sing q>2} Let $N\geq 1$ and $q>2$. If $u\in C^2(B_{r_0}\setminus\{0\})$ is a radial solution of $(\ref{I-1})$ in $B_{r_0}\setminus\{0\}$, the following dichotomy holds:\smallskip

\nind (i) either 
\bel{IV-4}
\lim_{r\to 0}r^qe^{u(r)}=Mq^q\quad\text{and }\,\lim_{r\to 0}ru_r(r)=-q,
\ee

\nind (ii) or there exists $u_0\in\BBR$ such that 
\bel{IV-5}
u(r)=u_0+c_{M,N,q}r^{\frac{q-2}{q-1}}(1+o(1))\quad\text{as }\,r\to 0,
\ee
where $c_{M,N,q}=\frac{q-1}{q-2}\left(\frac{(N-1)q-N}{M(q-1)}\right)^{\frac 1{q-1}}$ if $N>1$ and $c_{M,1,q}=-\frac{q-1}{q-2}\left(\frac{1}{M(q-1)}\right)^{\frac 1{q-1}}$,\smallskip

\nind (iii) or $u$ is regular at $0$ in the sense that there exists $u_0\in\BBR$ such that 
\bel{IV-7**}\dsps\lim_{r\to 0}u(r)=u_0\;\text{and }\;\dsps\lim_{r\to 0}u_r(r)=0,
\ee
\nind (iv) or $N=1$  there exists $u_0\in\BBR$ such that 
\bel{IV-7q>2}
\lim_{r\to 0}u(r)=u_0\in\BBR\quad\text{and }\lim_{r\to 0}u_r(r)=b\,\in\BBR_*.
\ee.
\es

For proving this result we need two intermediate lemmas.

\blemma{sing q>2-1} Assume $N\geq 1$ and $q>2$. If $u$ is a solution of $(\ref{II-1})$ in $(0,1)$, then \smallskip

\nind either for any $C\in (0,1)$ there exists $r_C>0$ such that\smallskip
\bel{IV-7***}r^qe^{u(r)}\geq CMq^q\;\text{ for all }r\geq r_C,
\ee
\nind or $u$ satisfies (ii), (iii) or (iv) in the previous statement.
\es
\Proof {\it Step 1: the case of decreasing solutions}. Assume $u_r<0$ near $0$. Then $u_r<0$ in $(0,1)$ by \rlemma{mono}. Following an idea of Serrin and Zou \cite{SeZo}, for any $C\in (0,1)$ we define the function $F_C$ by
\bel{IV-Fc}F_C(r)=e^u-CM|u_r|^q.
\ee
(i)- We first prove that given $C\in (0,1)$, $F_C$ has a constant sign near $0$. There holds
$$F'_{C}(r)=-e^u|u_r|+CMq|u_r|^{q-1}\left(M|u_r|^q+\frac{N-1}{r}|u_r|-e^u\right).
$$
At a point $r^*$ such that $F_{C}(r^*)=0$ we have
$$\BA{lll}\dsps
F'_{C}(r^*)=-CM|u_r|^{q+1}+CMq|u_r|^{q-1}\left(M(1-C)|u_r|^{q}+\frac{N-1}{r^*}|u_r|\right)\\[2mm]
\phantom{F'_{C}(r^*)}\dsps=-CM|u_r|^{q+1}+CMq|u_r|^{q}\left(M(1-C)|u_r|^{q-1}+\frac{N-1}{r^*}\right)\\[2mm]
\phantom{F'_{C}(r^*)}\dsps=CM|u_r|^{q}\left(qM(1-C)|u_r|^{q-1}+q\frac{N-1}{r^*}-|u_r|\right).
\EA$$
As a consequence, if furthermore $F'_C(r^*)<0$ we obtain that for $r=r^*$,
\bel{IV-7*}qM(1-C)|u_r|^{q-1}+q\frac{N-1}{r^*}<|u_r|.
\ee
Since $q>2$, this inequality implies that $|u_r(r^*)$ is bounded and more precisely $|u_r(r^*)|\leq A_{C,M,q}=(qM(1-C))^{\frac{1}{2-q}}$ which in turn implies 
$e^{u(r^*)}=CM|u_r(r^*)|^q\leq B_{C,M,q}=CM(qM(1-C))^{\frac{q}{2-q}}$. This implies also that $q\frac{N-1}{r^*}$ is bounded. Therefore, for 
$N>1$ this cannot happen and $F_C$ has a constant sign near $0$.\\
\nind Next, if $N=1$, $u_r$ is monotone near $0$ by \rlemma{mono}. If the function $F_C$ is oscillating near $0$, there exist two sequences $\{r_n\}$ and $\{\tilde r_n\}$ tending to $0$ such that $F_C(r_n)=F_C(\tilde r_n)=0$ with $F'_C(r_n)<0$ and $F'_C(\tilde r_n)>0$. This implies 
$|u_r(r_n)|<A_{C,M,q}$ and $|u_r(\tilde r_n)|>A_{C,M,q}$, which contradicts the monotonicity of $u_r$. As a consequence,  for $N=1$ also, $F_C$ keeps a constant sign near $0$.\smallskip

 \nind (ii)- Suppose that for any $C\in (0,1)$ there exists $r_C\in (0,1)$ such that $F_C(r)>0$ for $0<r<r_C$. Then $e^{\frac uq}+(CM)^\frac 1q u_r\geq 0$ on this interval which implies 
that $r\mapsto (CM)^\frac 1qr-qe^{-\frac {u(r)}q}$ in increasing on $(0,r_C)$. Since $u_r<0$, either $u(r)\to\infty$ when $r\to 0$, in which case  $(CM)^\frac 1qr-qe^{-\frac {u(r)}q}\geq 0$ and therefore 
$e^{u(r)}\geq CMq^qr^{-q}$ which is $(\ref{IV-7***})$. Or $u(r)\to u_0$ and the assumption $F_C>0$ implies  that $u_r$ is bounded on $(0,r_C]$. By standard regularity theory , $u_{rr}$ is also bounded. This implies that $u_r(r)$ has a limit when $r\to 0$, and this limit is necessarily zero from the equation if $N>1$. Therefore in that case $u$ satisfies 
$(\ref{IV-7q>2})$ or $(\ref{IV-7**})$.\smallskip

 \nind (iii)- Suppose now  that for some $C\in (0,1)$ there exists  $r_C\in (0,1)$ such that $F_C(r)<0$ for $0<r<r_C$. Then 
 $$-u_{rr}-\frac{N-1}{r}u_r+M(1-C)|u_r|^q\leq 0.
 $$
Set $\gm=M(1-C)$ and $W=-r^{N-1}u_r=r^{N-1}|u_r|$, then 
$$\gm r^{-(N-1)(q-1)}+W^{-q}W_r\leq 0.
$$
By integration, it implies that the function 
$$r\mapsto W^{1-q}(r)+\frac{\gm (q-1)}{(N-1)q-N}r^{N-(N-1)q}
$$
is nondecreasing. When $N>1$ we have that $(N-1)q-N>0$ since $q>2$ and we get a contradiction. When $N=1$, $W=-u_r$ and $r\mapsto |u_r(r)|^{1-q}-\gm(q-1)r$ is nondecreasing. Hence it admits  a limit $\ell\geq 0$ when $r\to 0$. If $\ell=0$ then $|u_r(r)|\to\infty$. Since $|u_r(r)|^{1-q}-\gm(q-1)r\to 0$, we have that $|u_r(r)|^{1-q}\geq \gm(q-1)r$. Hence 
$|u_r(r)|\leq (\gm(q-1)r)^{-\frac{1}{q-1}}$. By integration we obtain that $u(r)$ has a limit $u_0$. Since $-u_{rr}+M|u_r|^q(1+o(1))=0$ we deduce by integration that $(\ref{IV-5})$ holds. If $\ell>0$ then $|u_r(r)|\to \frac{\ell}{q-1}$ and $(\ref{IV-7q>2})$ holds.\smallskip

\nind {\it Step 2: the case of increasing solutions}. Since $u$ is bounded from above,  it follows by \rth {sup-2} that $|u_r(r)|\leq cr^{-\frac{1}{q-1}}$ near $0$. Since $q>2$ it follows that 
$u(r)$ admits a limit $u_0$ when $r\to 0$. The function $v=u-u_0$ satisfies
$$-v_{rr}-\frac{N-1}{r}v_r+M|v_r|^q=e^{u_0}(e^{v}-1):=\phi v
$$
where $\phi=e^{u_0}\frac{e^{v}-1}{v}\to e^{u_0}$ as $r\to 0$. Set $W=r^{N-1}v_r$, then
$$W_r=r^{N-1}\left(Mr^{-(N-1)q}W^q-\phi v\right).
$$
The function $u_r$ is monotone near $0$. Indeed, at a point $\tilde r$ where $u_{rr}(\tilde r)=0$ there holds
$$u_{rrr}(\tilde r)=\left(\frac{N-1}{\tilde r^2}-e^{u(\tilde r)}\right)u_r(\tilde r);
$$
as $u(r)\to u_0$ when $r\to 0$, $\frac{N-1}{ r^2}-e^u$ is positive near $0$ if $N\geq 2$, negative if $N=1$ and in any case this expression keeps a constant sign. \\
Either $u_r\to\infty$ when $r\to 0$ or $u_r$ has a finite limit. In the first case 
we have that $\phi v=o(v_r)$ when $r\to 0$. Thus we can write the equation satisfied by $v$ under the form
$$-r^{1-N}(r^{N-1}v_r)_r+M^*(r)|v_r|^q=0
$$
where $M^*(r)=M(1+o(1))$. This equation  is just a perturbation of  $-r^{1-N}(r^{N-1}v_r)_r+Mv_r|^q$ and it is explicitely integrable. If $N>1$ it yields $v_r(r)=\Gth_{M,N}r^{-\frac{1}{q-1}}(1+o(1))$ where 
$\Gth_{M,N}=\left(\frac{(N-1)q-N}{M(q-1)}\right)^{\frac 1{q-1}}$. We obtain $(\ref{IV-5})$. If $N=1$ we obtain a contradiction. In the case where $u_r$ has a finite limit, the function $u$ satisfies $(\ref{IV-7**})$ or $(\ref{IV-7q>2})$.\qeda\medskip

Next we give the precise behaviour of solutions of $(\ref{II-1})$ such that $r^qe^u$ is positively bounded from below. The result is obtained thanks to an energy function adapted from Leighton's method \cite{AnLei}.
\blemma{conver} Let $N\geq 1$ and $q>2$. If $u$ is a solution of $(\ref{II-1})$ such that $r^qe^{u}\geq C_1$ for some $C_1>0$ in a neighbourhood of $0$, then $(\ref{IV-4})$ holds.
\es
\Proof Since $u(r)\geq -q\ln r+c$ and $u_r$ has constant sign near $0$, it is negative. We set 
\bel{IV-8}X(t)=r^qe^{u(r)}\,,\; Y(t)=-r^{q+1}e^{u(r)}u_r(r)\quad\text{with }t=\ln r.
\ee
Then $(X,Y)$ satisfies 
\bel{IV-9+}\BA{lll}
X_t=qX-Y:=f(X,Y)\\[1mm]
Y_t=(q-N+2)Y-\frac{Y^2}{X}+e^{(2-q)t}\left(X^2-M\frac{Y^q}{X^{q-1}}\right)=g(X,Y,t).
\EA\ee
This system is a variant of system $(\ref{II-3})$ where the unknown are $X(t)$ and $\Phi(t)=-ru_r(r)$. Hence $Y(t)=X(t)\Phi(t)$.  Note that the use of $(X,Y)$ is the most natural transformation for the equation associated to $(\ref{III-5})$. We recall that Leighton's method relative to an autonomous system
$$\BA{lll}
\CX_t=\tilde f(\CX,\CY)\\[1mm]
\CY_t=\tilde g(\CX,\CY).
\EA$$
which has the property that the relation $\tilde f(\CX,\CY)=0$ is equivalent to $\CY=\tilde h(\CX)$ consists in analysing the variations of  the function 
$t\mapsto \CF(\CX(t),\CY(t))$ defined by
\bel{IV-10+}\CF(\CX,\CY)=\int_{h(\CX)}^{\CY}\tilde f(\CX,s)ds-\int_{h(\CX)}^{\CY}\tilde g(s,h(s))ds.
\ee
We adapt this methods to the non-autonomous system (\ref{IV-9+}) in considering first
$$\BA{lll}\dsps F(t)=\int_{qX(t)}^{Y(t)}f(X(t),s)ds-\int_{qX(t)}^{Y(t)}g(s,qs,t)ds\\[3mm]
\phantom{F(t)}\dsps =(N-2)q\frac {X(t)}2-\frac {X_t^2(t)}2+e^{(2-q))t}\left(Mq^q\frac{X^2(t)}{2}-\frac{X^3(t)}{3}\right).
\EA$$
In this expression the term $e^{(2-q))t}$ tends to $\infty$ when $t\to-\infty$. Therefore we replace $F$ by 
\bel{IV-G}G(t):=e^{(q-2)t}F(t)=Mq^q\frac{X^2(t)}{2}-\frac{X^3(t)}{3}+e^{(q-2))t}\left((N-2)q\frac {X(t)}2-\frac {X_t^2(t)}2\right).
\ee
 Then we obtain
 $$\BA{lll}\dsps G_t=(Mq^q-X)XX_t+(q-2)e^{(q-2))t}\left((N-2)q\frac {X}2-\frac {X_t^2}2\right)+e^{(q-2))t}X_t\left(\frac {(N-2)q}2-X_{tt}\right)\\[4mm]
 \phantom{G_t}\dsps =(Mq^q-X)XX_t+(q-2)e^{(q-2))t}\left((N-2)q\frac {X}2-\frac {X_t^2}2\right)\\[4mm]
 \phantom{G_t-------}\dsps+e^{(q-2))t}X_t\left(\frac {(N-2)q}2-(N-2)Y-\frac{X_t^2}{X}-e^{(2-q)t}\left(M\frac{Y^q}{X^{q-1}}-X^2\right)\right)
 \\[4mm]
 \phantom{G_t}\dsps=\left(Mq^qX-M\frac{Y^q}{X^{q-1}}\right)X_t+e^{(q-2))t}\Psi(t), \EA$$
where
  $$\BA{lll}\dsps
  \Psi(t)=(q-2)\left((N-2)q\frac X2-\frac{X_t^2}{2}\right)+X_t\left(\frac{(N-2)q}{2}-(N-2)Y-\frac{X_t^2}{X}\right)\\[4mm]
  \phantom{  \Psi(t)}\dsps
  =(q-2)\left((N-2)q\frac X2-\frac{X_t^2}{2}\right)+\frac{(N-2)q}{2}X_t-\frac{X_t^3}{X}+(N-2)(X_t-qX)X_t.
  \EA$$
Replacing $X_t$ by its value we obtain the following expression for $G_t(t)$
  \bel{IV-11*}
  G_t(t)=MX^{1-q}(qX-Y)(q^qX^q-Y^q)+e^{(q-2)t}\Psi(t).
\ee
Using the assumption and the bound from \rth{sup-1}, we have in a neighbourhood of $0$
$$C_1\leq r^qe^{u(r)}\leq C_2.
$$
This implies that $X(t)$ is bounded from above and from below. By the proof of \rlemma{sing q>2-1} we have for any $C<1$,
$CM|u_r|^q\leq e^{u(r)}$ near $0$, hence $|u_r|\leq \frac{C_3}{r}$ near $0$. Since $X(t)$  is bounded  when $t\to-\infty$ and $Y(t)=\phi(t)X(t)=-ru_r(r)X(t)$, we have also that $\Phi(t)$ and $Y(t)$ are bounded at $-\infty$. Hence $X_t(t)=qX(t)-Y(t)$ shares the same property. Furthermore
$\frac{X_t^3}{X(t)}=X^2(t)(q-\Phi(t))^2)$ is bounded, a fact which implies that the function $\Psi(t)$ is bounded. Noticing that  $(qX-Y)(((qX)^q-Y^q)\geq 0$, we obtain
$$G_t(t)=MX^{1-q}(qX-Y)(((qX)^q-Y^q)+e^{(q-2)t}\Psi\geq -Ce^{(q-2)t}\quad\text{for }t\leq 0.
$$
This implies that the function $t\mapsto G(t)+\frac{C}{q-2}e^{(q-2)t}$ is increasing. Hence either it tends to some finite $\ell$ or it tends to $-\infty$ when $t\to -\infty$. Since $X(t)$ is bounded from above and below and $Y(t)$ is bounded, we deduce that $\ell$ is finite. By the definition of $G(t)$ we have that $Mq^q\frac{X^2(t)}{2}-\frac{X^3(t)}{3}\to\ell$ when $t\to-\infty$, hence $X(t)$ converges to some 
$\gl$ such that $Mq^q\frac{\gl^2}{2}-\frac{\gl^3}{3}=\ell$, and $\ell\neq 0$ since $X(t)$ is bounded from below. Concerning $Y(t)$, either it has a limit $\Gl$ and $X_t(t)\to q\gl-\Gl$. But the only possible limit for $X_t(t)$ at $-\infty$ is zero, thus $\Gl=q\gl$. Or $Y(t)$ is not monotone, thus it oscillates at infinity and for a sequence $\{t_n\}$ of extremal points of $Y$ tending to $-\infty$,    we have that
$(q-N+2)Y(t_n)-\frac{Y^2(t_n)}{X(t_n)}=O(e^{(q-2)t_n})\to 0$  from the equation. Since $X(t_n)\to \gl$ and $e^{(q-2)t_n}\to 0$ it follows that $Y(t_n)$ admits a limit when $t_n\to\infty$, which implies in turn that 
 $Y^q(t_n)\to \frac{\gl^{q+1}}{M}$. Hence, even if it oscillates, $Y(t)$ admits a limit that we denote $\tilde \Gl$ at $-\infty$,  and $\tilde \Gl= \left(\frac{\gl^{q+1}}{M}\right)^{\frac 1q}$. Therefore 
 $\dsps\lim_{r\to 0}r^qe^{u(r)}=Mq^q$. Since $\Gf(t)=\frac{Y(t)}{X(t)}$, then  $\dsps\lim_{t\to-\infty}\Gf(t)=q=-\lim_{r\to 0}ru_r(r)$.\qeda\medskip
 
 In the next lemma we give precise estimates of the behaviour of such solutions by a method inspired from \cite{BiGaVe4}.
 
 \blemma{dev1} Let $N>1$ and $q>2$. If $u$ is a solution of $(\ref{II-1})$ satisfying $(\ref{IV-4})$, then
   \bel{IV-12}
u_{rr}(r)=\frac{q}{r^2}(1+o(1))\quad\text{as }r\to 0.
\ee
Furthermore
   \bel{IV-13}
u(r)=\ln\frac{Mq^q}{r^q}-\frac{N-2}{Mq^q(q-1)}r^{q-2}(1+o(1))\quad\text{as }r\to 0,
\ee
and
   \bel{IV-14}
u_r(r)=-\frac{q}{r}\left(1-\frac{(N-2)(q-2)}{Mq^q(q-1)}r^{q-2}(1+o(1))\right)\quad\text{as }r\to 0.
\ee
 \es
 \Proof In the proof we uses the system $(\ref{II-4})$ in $X(t)=r^qe^{u(r)}$ and $\Phi(t)=ru_r(r)$ with $t=\ln r$. We recall it below 
 $$\left\{\BA{lll}
 X_t=X(q-\Phi)\\[1mm]
 \Phi_t=(2-N)\Phi+e^{(2-q)t}(X-M|\Phi|^q).
\EA \right.$$
Then $\Phi>0$ since $u_r<0$. Note that $(Mq^q,q)$ is an equilibrium of $(\ref{II-4})$ only if $N=2$.\smallskip

\nind (i) {\it We claim that }
   \bel{IV-14*}
\lim_{t\to-\infty}e^{(2-q)t}(X(t)-M\Phi^q(t))=(N-2)q.
\ee
Set
$$\Psi(t)=e^{(2-q)t}(X(t)-M\Phi^q(t))=\Phi_t(t)+(N-2)\Phi(t).
$$
Then
$$\BA{lll}
\Psi_t=\Phi_{tt}+(N-2)\Phi_t=(2-q)e^{(2-q)t}(X-M\Phi^q)+e^{(2-q)t}(X(q-\Phi)-Mq\Phi^{q-1}\Phi_t)\\[1mm]
\phantom{\Psi_t}=(2-q)\Psi+e^{(2-q)t}X(q-\Phi)-Mqe^{(2-q)t}\Phi^{q-1}((2-N)\Phi+\Psi)\\[1mm]
\phantom{\Psi_t}=\left((2-q-Mqe^{(2-q)t}\Phi^{q-1}\right)\Psi+e^{(2-q)t}\left(X(q-\Phi)+Mq(N-2)\Phi^q\right)\\[1mm]
\phantom{\Psi_t}=e^{(2-q)t}\left(X(q-\Phi)+Mq(N-2)\Phi^q+\Psi\left((2-q)e^{(q-2)t}-Mq\Phi^{q-1}\right)\right).
\EA$$
If $\Psi$ is not monotone near infinity, at each $t^*<0$ where $\Psi_t(t^*)=0$ there holds
$$\left((q-2)e^{(q-2)t^*}+Mq\Phi^{q-1}(t^*)\right)\Psi(t^*)=X(t^*)(q-\Phi(t^*))+Mq(N-2)\Phi^q(t^*).
$$
Equivalently
$$\left(\frac{q-2}{Mq\Phi^{q-1}(t^*)}e^{(q-2)t^*}+1\right)\Psi(t^*)=\frac{X(t^*)(q-\Phi(t^*))}{Mq\Phi^{q-1}(t^*)}+(N-2)\Phi(t^*).
$$
By assumption $\Phi(t)\to q$ and $X(t)\to Mq^q$ as $t\to-\infty$. If $t^*=t_n\to-\infty$, we have that  $\Psi(t_n)\to (N-2)q$, then $\Psi(t)\to (N-2)q$. \\
If $\Psi$ is monotone near infinity, then it admits a limit $L\in [-\infty,\infty]$. Since $\Phi_t(t)\to L-(N-2)q$ and $\Phi(t)\to q$ when $t\to-\infty$, then necessarily $L$ is finite, $\Phi_t(t)\to 0$  and finally   $L=(N-2)q$ as in the previous case. It follows that 
$$\lim_{t\to-\infty}\Psi(t)=(N-2)q=\lim_{t\to-\infty}e^{(2-q)t}(X(t)-M\Phi(t)),
$$
which can be written
$$X(t)-M\Phi(t)=(N-2)qe^{(q-2)t}(1+o(1))\quad\text{as }t\to-\infty,
$$
which is the claim.\smallskip

\nind (ii) {\it End of the proof. } Relation (\ref{IV-14*}) can be expressed by
\bel{IV-14**}r^2(e^{u(r)}-M|u_r(r)|^q)=(N-2)q(1+o(1))\quad\text{as }r\to 0.
\ee
By  $(\ref{II-1})$ we obtain $(\ref{IV-12})$ since
$$\BA{lll}\dsps
u_{rr}=\frac{N-1}{r}|u_r|-e^u+M|u_r|^q\\[4mm]
\phantom{u_{rr}}
\dsps =\frac{(N-1)q}{r^2}(1+o(1))-\frac{(N-2)q}{r^2}(1+o(1))=\frac{q}{r^2}(1+o(1)).
\EA$$
Since $M|u_r(r)|^q=e^{u(r)}(1+o(1))$ when $r\to 0$, we define $W$ by $u_r(r)=-\frac{e^{\frac{u(r)}{q}}}{M^{\frac 1q}}(1+W(r))$. Then $W(r)\to 0$ and we have when $r\to 0$,
$$r^2(e^{u(r)}-M|u_r(r)|^q)=r^2e^{u(r)}(1-(1+W(r))^q)=-qr^2e^{u(r)}W(r)(1+o(1)).
$$
This implies
$$W(r)=-\frac{(N-2)e^{-u(r)}}{r^2}(1+o(1))=-\frac{(N-2)r^{q-2}}{Mq^q}(1+o(1)).
$$
Therefore
   \bel{IV-15}-u_r(r)=\frac{e^{\frac{u(r)}{q}}}{M^\frac 1q}\left(1-\frac{(N-2)r^{q-2}}{Mq^q}(1+o(1))\right),
\ee
which can be written as
$$M^\frac 1q e^{-\frac{u(r)}{q}}u_r(r)+1-\frac{(N-2)r^{q-2}}{Mq^q}(1+o(1))=0.
$$
By integration, 
$$-qM^\frac 1q e^{-\frac{u(r)}{q}}+r\left(1-\frac{(N-2)r^{q-2}}{M(q-1)q^q}(1+o(1))\right)=0.
$$
This yields the expansion of $e^{u(r)}$,
   \bel{IV-16}
e^{u(r)}=\frac{Mq^q}{r^q}\left(1+\frac{(N-2)r^{q-2}}{M(q-1)q^{q-1}}(1+o(1))\right),
\ee
from which follows 
   \bel{IV-17}\BA{lll}\dsps
u(r)=\ln\frac {Mq^q}{r^q}+\frac{(N-2)r^{q-2}}{Mq^q(q-1)}(1+o(1))\quad\text{as }r\to 0.
\EA\ee
Combining $(\ref{IV-15})$ and $(\ref{IV-16})$ we obtain
   \bel{IV-18}\BA{lll}\dsps
-u_r(r)=\frac{q}{r}\left(1-\frac{(N-2)(q-2)r^{q-2}}{Mq^q(q-1)}(1+o(1))\right)\quad\text{as }r\to 0.
\EA\ee
This ends the proof.\qeda\medskip

\nind\Remark By $(\ref{IV-16})$ we obtain the 
   \bel{IV-19}\BA{lll}\dsps
r^qe^{u(r)}>Mq^q \, \left(\text{resp }\; r^qe^{u(r)}<Mq^q\right) \,\text{ if }\;N\geq 2\,\left( \text{resp }\;N=1\right).
\EA\ee
\subsection{Existence of singular solutions of eikonal type}
In this section we prove existence results for eikonal type solutions. The proof is difficult and in dimension $N\geq 2$ it is based upon the construction for $N=1$.
\subsubsection{The case $N=1$}
\bth{uniq-ei} Let $N= 1$ and $q>2$. Then there exists one and only one solution $u^*$ of $(\ref{II-1})$ in $(0,\infty)$ such that $(\ref{IV-4})$ holds. Furthermore $\dsps u^*=\lim_{n\to\infty}u_n$ where $u_n$ is the regular solution of $(\ref{II-1})$ in $(0,\infty)$ such that $u_{n\,r}(0)=0$ and $u_n(0)=n$.
\es
\Proof {\it 1- Uniqueness}. We already know that such a solution $u$ is decreasing and satisfies $\dsps \lim_{r\to\infty}u(r)=-\infty$. Hence $r\mapsto u(r)$ is a decreasing diffeomorphism from $(0,\infty)$ onto $(-\infty,\infty)$;  we consider $r$ as a function of $u$ and set 
$$z(u)=u_r^2(r)=u_r^2(r(u)).$$ 
Then $z$ is defined on 
$\BBR$ and there holds
$$\frac {dz}{du}=2u_ru_{rr}\frac{dr}{du}=2u_{rr},
$$	
hence
   \bel{IV-21}\BA{lll}\dsps
\frac{dz}{du}=2Mz^{\frac q2}-2e^u.
\EA\ee
The associated system in $(r,z)$ as functions of $u$ is
   \bel{IV-22*}\left\{\BA{lll}\dsps
z_u=2Mz^{\frac q2}-2e^{u}\\[0mm]\dsps
r_u=-\frac{1}{\sqrt z}.
\EA\right.\ee
By $(\ref{IV-4})$  we have that
   \bel{IV-22}z(u)=\frac{e^{\frac {2u}{q}}}{M^{\frac 2q}}(1+o(1))\quad\text{and }\;r(u)=M^{\frac 1q}qe^{-\frac uq}(1+o(1))\quad\text{as }\;u\to\infty.
\ee
The point is that there exists a {\it unique} solution of $(\ref{IV-21})$ satisfying $(\ref{IV-22})$. To see that, 
let $z_1$ and $z_2$ be two such solutions,  then by substracting the two corresponding equations in $z_j$ (j=1,2), the term $-2e^u$ disappears and we obtain 
$$\BA{lll}\dsps\frac {d(z_1-z_2)}{du}=2M\left(z_1^{\frac q2}-z_2^{\frac q2}\right)=Mq\xi^{\frac {q-2}2}(z_1-z_2),
\EA$$
where $\xi:=\xi(u)=\gth z_1(u)+(1-\gth)z_2(u)$ for some $\gth\in [0,1]$. Because of  $(\ref{IV-4})$ we have that 
$$Mq\xi^{\frac {q-2}2}=qM^\frac 2qe^{(1-\frac 2q)u}(1+o(1))\quad\text{as }u\to\infty.$$
Hence the function
$$u\mapsto H(u,v):=e^{-Mq\int_{v}^{u}\xi(s)^{\frac {q-2}2}ds}(z_1(u)-z_2(u))
$$
is constant for any $v$. Since by l'Hospital's rule
$$-Mq\int_{v}^{u}\xi(s)^{\frac {q-2}2}ds=-\frac{Mq^2}{q-2}e^{(1-\frac 2q)u}(1+o(1)))\quad\text{as }u\to\infty.
$$
it follows from $(\ref{IV-22})$ that $H(u,v)\to 0$ as $u\to\infty$, hence $H(v,v)=0$ which implies $z_1(v)=z_2(v)$. Therefore if $u_1$ and $u_2$ are solutions, the uniqueness 
of $z_1$ and $z_2$ implies $u_{1\,r}(r)=u_{2\,r}(r)<0$ for any $r>0$ and $u_{1}(r)=u_{2}(r)+c$. The fact that $u_1$ and $u_2$ satisfy $(\ref{II-1})$ implies $c=0$.\smallskip
\smallskip

\nind {\it 2- Existence}. Equation $(\ref{II-1})$ reduces to 
   \bel{IV-23}\BA{lll}\dsps
-u_{rr}+M|u_r|^q-e^u=0
\EA\ee
equivalently to the {\it autonomous }system of order 2,
   \bel{IV-24}\left\{\BA{lll}\dsps
u_{r}=-v\\[0mm]\dsps
v_{r}=e^u-M|v|^q.
\EA\right.\ee
Let $a_1\geq 0$ and denote by $u_a$ any regular solution of $(\ref{IV-23})$ such that $u(0)=a>a_1$ and $u_r(0)=0$. Since $u_a$ is monotone it is decreasing and necessarily it tends to $-\infty$ at infinity; from \rprop{glob} there exists a unique $\gr_1$ such that $u_a(\gr_1)=0$. \smallskip

\nind{\it 2- Step 1. We claim that there exists $K=K_{a_1,M,q}>0$ such that }
   \bel{IV-25}\BA{lll}\dsps
|u_{a\,r}(\gr_1)|\leq K.
\EA\ee
We consider the function  $F_C$ defined in $(\ref{IV-Fc})$ for $C=\frac 12$ and $F_{\frac 12}(r)=e^{u_a(r)}-\frac 12M|u_{a\,r}(r)|^q$. Then $F_{\frac 12}(0)=e^a$ and 
$$F_{\frac 12}(\gr_1)=e^{a_1}-\tfrac 12M|u_{a\,r}(\gr_1)|^q.$$
If $F_{\frac 12}(\gr_1)\geq 0$ we obtain 
$$|u_{a\,r}(\gr_1)|\leq \left(\frac{2e^{a_1}}{M}\right)^{\frac 1q}.
$$
If $F_{\frac 12}(\gr_1)< 0$ there exists $r_0\in (0,\gr_1)$ such that $F_{\frac 12}(r_0)=0$ and $F'_{\frac 12}(r_0)<0$ and by $(\ref{IV-7*})$, we deduce $|u_{a\,r}(r_0)|^q>0$, then
$2^{-1}qM|u_{a\,r}(r_0)|^{q-1}\leq |u_{a\,r}(r_0)|$, hence 
$$|u_{a\,r}(r_0)| \leq (2^{-1}qM)^{\frac{1}{q-2}}:=\tilde K_{a_1,M,q},$$
and
$$-u_{a\,rr}(r_0)=e^{u_a(r_0)}-M|u_{a\,r}(r_0)|^q=-2^{-1}M|u_{a\,r}(r_0)|^q<0,
$$
Since $u_{a\,rr}(0)=-e^{a}<0$, there exists $r_2\in (0,r_0)$ such that $u_{a\,rr}(r_2)=0$, and this point is unique since $-u_{a\,rrr}(r_2)=e^{u_a(r_2)}u_{a\,r}(r_2)<0$. Then $u_{a\,rr}>0$ on 
$(r_2,\gr_1)$. This implies that $u_{a\,r}$ is increasing on $(r_2,\gr_1)$ and thus $u_{a\,r}(r_2)<u_{a\,r}(\gr_1)$. Consequently 
$$|u_{a\,r}(\gr_1)|\leq |u_{a\,r}(r_0)|\leq \tilde K_{a_1,M,q}
$$
Estimate $(\ref{IV-25})$ follows with $K=\max\{\left(\frac{2e^{a_1}}{M}\right)^{\frac 1q},\tilde K_{a_1,M,q},\}$.
\smallskip

\nind {\it 2- Step 2: We claim that there exist singular solutions.} We denote by $\CH$ the vector field 
$$\CH(u,v)=(-v,e^u-M|v|^q):=(\CH_1(u,v),\CH_2(u,v))
$$
in the phase plane $Q:=\{(u,v):v>0\}$. we denote by $\CC$ the curve $\{(u,v)\in Q: Mv^q=e^u\}$ and the two regions
$$\CR:=\{(u,v)\in Q: Mv^q>e^u\}\,\text{ and }\;\CS:=\{(u,v)\in Q: Mv^q<e^u\}.
$$
On $\CC$ the vector field is horizontal and entering in $\CR$. Thus this region is positively invariant. On the axis $v=0$ the vector field is vertical and inward to $Q$.\smallskip

\nind {\bf-} A regular trajectory $\CT_a:=\{(u_a,-u_{a\,r})\}_{r>0}$ starts from $(a,0)$ with a vertical slope. Since there exists a unique $r_2>0$ such that $u_{a\,rr}(r_2)=0$, 
$\CT_a$ intersects $\CC$ in $Q$ and remains in $\CR$ which is indeed positively invariant. Conversely, any trajectory issued from a point 
$(\bar u,\bar v)\in\CC$ for $r=\bar r$ is included in  $\CS$ for $r<\bar r$ as long as it remains in $Q$. Since $\CH$ admits no equilibrium in $\CS$, this backward trajectory 
cannot remain bounded for $r<\bar r$. If this backward trajectory does not intersect the axis $v=0$, then $u(r)\to\infty$ when $r<\bar r$ decreases to the infimum of the maximal interval 
of existence and $v(r)$ admits a finite limit. From $(\ref{IV-24})$ this implies that $v_r(r)\to \infty$, contradiction. Then the trajectory intersect $v=0$ at some $\tilde r<\bar r$. This implies that $\bar u(r-\tilde r)$ is the regular trajectory issued from 
$\bar u(\tilde r)$.\smallskip

\nind {\bf-}  We first take $a_1=0$. Hence any regular trajectory cuts the axis $u=0$ at a point $c$ such that $0<c\leq K_{0,M,q}$. Therefore, any trajectory through a point $(0,c)$ with 
$c>K_{0,M,q}$ is not a regular one. We denote by $Reg$ the set of $c>0$  such that the trajectory through $(0,c)$ is regular and put
\bel{IV-25*}
c^*=\sup\{c>0: c\in Reg\}.
\ee
By the implicit function theorem, the set of $Reg$  is open and it is an interval since two trajectories cannot intersect. Therefore $c*$ is not in this set. Let $\CT_*$ be the backward trajectory through 
$(0,c^*)$. It is does not intersect $\CC$, thus it is located in the region $\CR$ where $Mv^q>e^u$. This means $u_re^{-\frac uq}+M^{-\frac 1q}< 0$, thus the function 
$r\mapsto M^{-\frac 1q}r-qe^{-\frac{u(r)}{q}}$ is nonincreasing, consequently, for such a trajectory, the variable $r$ is bounded from below. Therefore any solution $u$ with trajectory $\CT^*$ is defined on a maximal interval $(R^*,\infty)$ with $R^*>-\infty$. Now the function $r\mapsto u^*:=u(r-R^*)$ which a solution is defined on the maximal interval $(0,\infty)$, it is nonincreasing, thus it is necessarily singular and the trajectory $\CT_a$ of any regular solution lies below $\CT^*$ in the half-plane $Q$. By \rth{sing q>2} we have two possibilities:\smallskip

\nind {\it (i)} Either there exists $a_1>0$ such that $u^*(r)\to a_1$ and $u^*_r(r)=-\left(\frac{1}{M(q-1)}\right)^{\frac1{q-1}}r^{-\frac 1{q-1}}(1+o(1))$ when $r\to 0$. In such a case, for any $a>a_1$, the regular solution $u_a$ which satisfies 
$u_a(\gr_1)=a_1$ for some $\gr_1>0$ satisfies also  $|u_{a\,r}(\gr_1)|\leq K_{a_1,M,q}$  by $(\ref{IV-25})$. This means that in the phase plane, on the trajectory $\CT_a$, there holds $v_a\leq K_{a_1,M,q}$ at the point 
$(a_1,v_a)=(u_a(\gr_1),-u_{r\,a}(\gr_1))$. Consider now  the trajectory passing through $(a_1,1+K_{a_1,M,q})$, it is below $\CT^*$ and above $\CT_a$ for all $a>a_1$. Therefore it intersects the axis 
$u=0$ at some $\tilde c<c^*$ but which is also larger or equal to $\sup \left\{c:c\in Reg\right\}$ which is $c^*$. This is a contradiction.\smallskip

\nind {\it (ii)} Or $\dsps \lim_{r\to 0}u^*(r)=\infty$. By \rth{sing q>2} this implies that $(\ref{IV-4})$ holds.\smallskip

\nind {\it 2-Step 3: Convergence of the regular solutions.} We consider the regular solutions $u_n$, this means $u_n(0)=n$ and $u_{n\,r}(0)=0$. Let $b\leq 0$, since $u_n$ is decreasing and tends to $-\infty$ when $r\to\infty$ there exists $\gr_{n,b}>0$ such that 
$u_n(r)>b$ on $[0,\gr_{n,b})$ and $u_n(\gr_{n,b})=0$. Similarly $u^*(r)>b$ on $[0,\gr_{*,b})$. We can write
$$\gr_{n,b}=\int_b^n\frac{du}{v_n(u)}=\int_b^n\frac{du}{\chi_{(0,n)}(u)v_n(u)}\;\text{ and }\gr_{*,b}=\int_b^\infty\frac{du}{v^*(u)}.
$$
Because two different trajectories cannot intersect, $v_n(u)$ is an increasing function of $n$. Hence, by the monotone function theorem 
$$\gr_{*,b}=\inf_{n}\gr_{n,b}=\lim_{n\to\infty}\gr_{n,b}.
$$
This implies that all the regular solutions are greater than $b$ on $[0,\gr_{*,b})$. By \rth{sup-1}-(i) we have for any $R>0$
$$e^{u_n(r)}\leq Cr^{-q}\quad\text{for all }r\in (0,R]
$$
where $C>0$ depends on $R$ but not on $n$ since $u_n$ is decreasing on $[0,\infty)$, therefore there also holds 
$$0\leq e^{u_n-b}\leq C e^{-b}r^{-q}\quad\text{on }\,[0,\gr_{*,b}].
$$
From $(\ref{III-7})$ in \rth{sup-2} we have
$$|u_{n\,r}(r)|\leq C\left(r^{-\frac 1{q-1}}+r^{-1}+1\right)\quad\text{for all }r\in (0,\tfrac R2],
$$
where $C>0$ is independent of $n$ and thus
$$|u_{n\,rr}(r)|\leq C\left(r^{-q}+r^{-\frac{q-2}{q-1}}+r^{-2}+1\right)\quad\text{on }\,[0,\gr_{*,b}].
$$
There exists a subsequence $\{n_j\}$ and a solution $U_b\geq b$ of $(\ref{IV-23})$ such that $\{u_{n_j}\}$ converges to $U_b$ in the $C^1_{loc}((0,\gr_{*,b}))$-topology. The corresponding trajectory 
$\CT_{U_b}$ is below $\CT_*$, thus $\CT_{U_b}$ crosses the axis $u=0$ at a point which is necessarily $c^*$ from $(\ref{IV-25*})$, hence $\CT_{U_b}=\CT_*$. Because of uniqueness, see Step-1, the whole sequence 
$\{u_{n_j}\}$ converges to $u^*$.\qeda
\subsubsection{The case $N\geq 2$}
\bth{uniq-ei-2} Let $N> 1$ and $q>2$. Then there exists at least one solution $u^*$ of $(\ref{II-1})$ in $(0,\infty)$ such that $(\ref{IV-4})$ holds. 
\es
\Proof Since the proof uses the result in dimension $1$, we will denote by $u^{(N)}$ or $u^{(1)}$ the solutions in N-dim or in 1-dim. Thus $u^{(N)}_a$ and $u^{(1)}_a$ denote the regular solutions respectively in $\BBR^N$ and in $\BBR$ with initial data $a$. For $b\leq 0$ we denote by $[0,\gr^{(N)}_{a,b})$ the maximal interval where $u^{(N)}_a>b$.\smallskip

\nind{\it Step 1: We claim that $\gr^{(N)}_{a,b}>\gr^{(1)}_{a,b}$ and $u^{(N)}_a>u^{(1)}_a$ on $(0,\gr^{(1)}_{a,b})$ for $N>1$.} Since $r\mapsto u^{(j)}(r)$ ($j=1,N$) is decreasing, we set $z^{(j)}(u)=\left(u^{(j)}_r\right)^2\left((u^{(j)})^{-1}(u)\right)$. Then 
$$\frac{dz^{(N)}}{du}-2M(z^{(N)})^{\frac q2}+2e^{u}=\frac{N-1}{r^{(N)}}\sqrt{z^{(N)}}
$$
for $N>1$, and
$$\frac{dz^{(1)}}{du}-2M(z^{(1)})^{\frac q2}+2e^{u}=0,
$$
we obtain that 
\bel{IV-26}
\frac {z^{(1)}-z^{(N)}}{du}=2M\left((z^{(1)})^{\frac q2}-(z^{(N)})^{\frac q2}\right)-\frac{N-1}{r^{(N)}}\sqrt{z^{(N)}}<2M\left((z^{(1)})^{\frac q2}-(z^{(N)})^{\frac q2}\right).
\ee
we have also
\bel{IV-27}
u^{(j)}_a(r)=a+\int_0^rs^{1-j}\int_0^s\left(M|u^{(j)}_{a\,r}|^q-e^{u^{(j)}_a}\right) t^{j-1}dtds,
\ee
hence
$$u^{(j)}_a(r)=-\frac{re^{a}}{j}(1+o(1))\quad\text{as }r\to 0.
$$
By integration we obtain in particular $u^{(N)}_a(r)>u^{(1)}_a(r)$ near $r=0$, and
$$z^{(1)}(u)-z^{(N)}(u)=\frac{e^{2a}(N^2-1)}{N^2}(1+o(1))\quad\text{as }u\to a.$$
Let $\hat a$ be the infimum of the $u\in (b,a)$ such that $z^{(1)}(u)-z^{(N)}(u)>0$.
If $\hat a>b$, then $z^{(1)}(\hat a)-z^{(N)}(\hat a)=0$ and $
\frac{d(z^{(1)}-z^{(N)})}{dr}(\hat a)\geq 0$, which contradicts $(\ref{IV-26})$, hence $z^{(1)}(u)>z^{(N)}(u)$ on $(b,a)$. By assumption $u^{(N)}(\gr^{N}_{a,b})=b$. This implies that 
$(u^{(N)}_a-u^{(1)}_a)_r>0$ on the interval $(0,\min\{\gr^{(N)}_{a,b},\gr^{(1)}_{a,b}\})$. By integration $u^{(N)}_a>u^{(1)}_a$ on this interval, which implies $\gr^{(N)}_{a,b}>\gr^{(1)}_{a,b}$ and 
$u^{(N)}_a>u^{(1)}_a$ on $(0,\gr^{(1)}_{a,b}]$, which is the claim.\smallskip

\nind{\it Step 2: Convergence of the regular solutions.}  Because the upper estimates of \rth{sup-1}-(i) and \rth{sup-2} hold independently of $n$, one can extract a subsequence $\{n_j\}$ such that $\{u_{n_j}^{N}\}$ converges in the $C^1_{loc}((0,\gr^{(1)}_{*,b}))$-topology to a function $U_{b}^{*(N)}>b$ which is a solution of $(\ref{II-1})$ on $(0,\gr^{(1)}_{*,b})$, is larger than the 1-dimensional solution 
$U_{b}^{*(1)}>b$. Thus it is singular and by \rth{sing q>2} it satisfies 
$$\lim_{r\to 0}r^qe^{U_{b}^{*(N)}(r)}=Mq^q\quad\text{and }\;\lim_{r\to 0}rU_{b\,r}^{*(N)}(r)=-q.
$$
This solution can be extended as a global solution on $(0,\infty)$ by \rprop{glob}.This ends the proof.\qeda\medskip

\nind\Remark The uniqueness of the singular solution of eikonal type when $N\geq 2$ is a challenging question. The remarkable fact is that all the solutions of this type have the same asymptotic expansion up to any order. This result is proved in Appendix.\medskip

\nind\Remark The Step 1 of our proof is an adaptation to equation $(\ref{II-1})$ of a method introduced by Voirol  \cite [Proposition 1.7]{Vo}
dealing with the Chipot-Weissler eqquation
\bel{IV-CW}
-u_{rr}-\frac{N-1}{r}u_r+|u_r|^q-u^p=0.
\ee

\subsection{Existence of singular solutions of Hamilton-Jacobi type}

 In this section we prove the existence of singular solutions with the same behaviour as the one of the solutions of the Hamilton-Jacobi equation.

\bth{q>2-beta} Let $N= 1$ and $q>2$. Then for any $u_0\in\BBR$  there exists at least one solution $u$ of $(\ref{II-1})$ in $(0,\infty)$ satisfying
\bel{IV-28*}
u(r)=u_0-\frac{q-1}{q-2}\left(\frac{1}{M(q-1)}\right)^{\frac 1{q-1}}r^{\frac{q-2}{q-1}}(1+o(1))\quad\text{when }\;r\to 0.
\ee
Furthermore $u$ is decreasing on $(0,\infty)$.
\es
\Proof We still use system $(\ref{IV-24})$. For any $w\in\BBR$, we denote by $(w,c^*_w)$ the intersection of $\CT_*$ with the axis $u=w$. For $k_0>c^*_w$ the trajectory $\CT_{k_0,w}$ going through $(w,k_0)$ is singular and differs from $\CT_*$. By \rth{sing q>2} $\CT_{k_0,w}$ has a vertical asymptote $u=u_0$ for some $u_0>w$. Because two trajectories cannot intersect, the mapping $k_0\mapsto f(k_0)=u_0$ is decreasing from $(c_w^*,\infty)$ to $(w,\infty)$. It is certainly continuous because of the non-intersection condition and onto from $(c_w^*,\infty)$ to $(w,\infty)$. This implies the claim for any $u_0>w$ and for any $w\in\BBR$.
 \qeda\medskip
 


 \bth{q>2-beta*} Let $N\geq 2$ and $q>2$. Then for any $u_0\in\BBR$ there exists at least one solution $u$ of $(\ref{II-1})$ in $(0,\infty)$ satisfying
\bel{IV-28}
u(r)=u_0+\frac{q-1}{q-2}\left(\frac{N(q-1)-N}{M(q-1)}\right)^{\frac 1{q-1}}r^{\frac{q-2}{q-1}}(1+o(1))\quad\text{when }\;r\to 0.
\ee
The function $u$ is increasing near $0$.\es
\Proof We use the system $(\ref{II-12})$ with $Z=-\frac{re^u}{u_r}$, $V=r|u_r|^{q-1}$ and $\Gf=-r u_r$. The solution $u$ we look for satisfies $\dsps \lim_{r\to 0}u(r)=u_0$ and $\dsps \lim_{r\to 0}r^{\frac 1{q-1}}u_r(r)=k$ for some real number $k$. Since we search increasing solutions, $k\geq 0$, $\Phi$ and $Z$ are non-positive and the system becomes
\bel{IV-29**}\left\{\BA{lll}
Z_t=Z(N-\Phi-MV-Z)\\
\dsps
V_t=V\left(N-(N-1)q+(q-1)(Z+MV)\right)\\
\Phi_t=\Phi\left(2-N+Z+MV\right).
\EA\right.\ee

 \nind{\it Step 1: The linearised problem}. In this system the relation $V=e^{(2-q)t}|\Gf|^{q-1}$ holds and thus
$$Z=e^{(2-q)t}\frac{X}{\Gf}\,,\; X=\frac{|\Gf|^{q-1}\Gf}{V}\,\text{ and }\;  e^{qt}X=ZV^{\frac{2}{q-2}}|\Phi|^{-\frac q{q-2}}.
$$
The system $(\ref{IV-29**})$ admits
$P_0=(0,V_0,0)$ for stationary point with $V_0=\frac{(N-1)q-N}{M(q-1)}$. Setting
$V=V_0+\bar V$ the linearised system is
\bel{IV-29}\left\{\BA{lll}\dsps
Z_t=\frac q{q-1}Z\\[2mm]
\dsps \bar V_t=(q-1)V_0(M\bar v+Z)\\[1mm]
\dsps \Phi_t=\frac{q-2}{q-1}\Phi.
\EA\right.\ee
The eigenvalues of the linearised system $(Z,\bar V,\Phi)$ at $(0,0,0)$ are $(\frac{q}{q-1},(N-1)q-N,\frac{q-2}{q-1})$. They are all positive. Hence there exists a 
neighbourhood  $\CV$ of $P_0$ such that all trajectories with a point in $\CV$ converge to $P_0$ when $t\to-\infty$.\smallskip

\nind {\it Step 2: End of the proof}. We recall that $u$ is a solution of $(\ref{II-1})$ if and only if $(X,\Gf)$ defined by $(\ref{II-2})$ satisfies the system $(\ref{II-4})$ but this system is not equivalent to  $(\ref{II-12})$. Let $\gs>0$ such that 
the ball $B_\gs(P_0)$ is in the attractive basin of $P_0$ when $t\to-\infty$. If $P^*:=(Z^*,V^*,\Phi^*) \in B_\gs(P_0)$ with $Z^*,\Phi^*<0$ we denote by $\CT_{P^*}$ the backward trajectory of 
$P^*$ and we set 
$$H(t)=\frac{V^{\frac 1{q-1}}}{\Phi}e^{\frac{q-2}{q-1}t}.$$
Then
$$\BA{lll}\dsps
(q-1)\frac {H_t}{H}=q-2+\frac{V_t}{V}-(q-1)\frac{\Phi_t}{\Phi}\\[2mm]
\phantom{(q-1)\frac {H_t}{H}}\dsps =q-2+N-q(N-1)+q-1)(MV+Z)-(q-1)(2-N+MV+Z)\\[2mm]
\phantom{(q-1)\frac {H_t}{H}}\dsps = q-2+N-(N-1)q+(q-1)(N-2)=0.
\EA$$
Hence the function $H$ is constant. This implies that $V=be^{(2-q)t}|\Phi|^{q-1}$ for some $b>0$. It is easy to check by computation that $X=\frac{Z|\Phi|^{q-1}\Phi}{V}=b^{-1}e^{(q-2)t}Z\Phi$ satisfies  the following system,
$$\left\{\BA{lll}\dsps
X_t=(q-\Phi)X\\[2mm]
\Phi_t=-(N-2)\Phi+be^{(2-q)t}(X-M|\Phi|^q).
\EA\right.$$
We define $a=\frac1{2-q}\ln b$, $\gt=t+a$, $X^{(a)}(\gt)=X(\gt-a)$ and $\Phi^{(a)}(\gt)=\Gf(\gt-a)$. Then $\Phi^{(a)}(\gt)\sim -V_0^{\frac 1{q-1}}e^{\frac{q-2}{q-1}t}$ as $\gt\to-\infty$ and 
$(X^{(a)},\Phi^{(a)})$ satisfies the system
$$\left\{\BA{lll}\dsps
X^{(a)}_t=(q-\Phi^{(a)})X^{(a)}\\[2mm]
\Phi^{(a)}_t=-(N-2)\Phi^{(a)}+e^{(2-q)\gt}(X^{(a)}-M|\Phi^{(a)}|^q).
\EA\right.$$
By \rlemma{L-1}, equivalently the function $\gr\mapsto u^{(a)}(\gr)=\ln(\gr^{-q}X^{(a)})(\ln\gr)$ satisfies $(\ref{II-1})$, with $\gr u^{(a)}_\gr=-\Phi^{(a)}(\gt)$ and 
$\gr (u^{(a)}_\gr)^{q-1}=e^{(2-q)\gt}(-\Phi^{(a)})^{q-1}$, thus $\dsps\lim_{\gr\to 0}\gr (u^{(a)}_\gr(\gr))^{q-1}=V_0$. Since $q>2$ the function $u^{(a)}_\gr$ is integrable near 
$0$ and there exists some $u_0$ such that $\dsps\lim_{\gr\to 0}u^{(a)}(\gr)=u_0$. This implies that 
$$u(r)-u_0=\frac {q-1}{q-2}V_0^{\frac 1{q-1}}r^{\frac{q-2}{q-1}}(1+o(1))\quad\text{as }r\to 0.
$$
Since $D\CG(P_0)$ is invertible the flows of $\CG$ is conjugate to the one of $D\CG(P_0)$ in a neighbourhood of $P_0$. hence there exists a $C^2$ diffeomorphism $\Gth$ from 
$B_\gs(P_0)$ (up to reducing $\gs$ in order $B_\gs(P_0)$ is a subset of the basin of attraction of $P_0$) to a neighbourhood $\CV$ of $0$ such that $\CG-\CG(P_0)=\CV\circ D\CF(P_0)$, and the trajectories of $(\ref{IV-29**})$ in $B_a(P_0)$ (for some $a>0$) are in one to one correspondence 
via $\CV$ with the trajectories of 
\bel{IV-31}(Z_t,V_t,\Phi_t)=D\CF(P_0)((Z,V,\Phi))
\ee
Since all the trajectories of $(\ref{IV-31})$ converge to $0$ when $t\to-\infty$, all the trajectories of $(\ref{IV-29**})$ issued from $B_\gs(P_0)$ converge to $P_0$ when $t\to-\infty$. This ends to proof.
Given arbitrary  coefficients $C_j$ ($j=1,2,3$) with $C_3>0$, the solutions of the linearized system $(\ref{IV-31})$ are expressed by
$$\left\{\BA{lll}Z(t)=C_1e^{\frac {qt}{q-1}}\\
\dsps \bar V(t)=C_2e^{(N(q-1)-N)t}+\frac{(q-2)C_1}{M(q-1)}e^{\frac{qt}{q-1}}\\
\Phi(t)=C_3e^{\frac {(q-2)t}{q-1}},
\EA\right.
$$
in general, with a standard modification if $N=3$ and $q=\frac{3+\sqrt 3}{2}$ or $N=2$ and $q=2+\sqrt 2$ since in these two cases $\frac{q}{q-1}=(N-1)q-N$. If $(V,\Gf)$ satisfies  $(\ref{II-12})$, then the equivalence of trajectories yields
$$V(t)=be^{(2-q)t}\Gf^{q-1}(t)\sim be^{(2-q)t}C_3^{q-1}e^{(q-2)t}=bC_3^{q-1},
$$
and clearly $b=V_0C_3^{1-q}=e^{(2-q)a}$ with the previous notations. Hence $e^{qa}=(V_0C_3^{1-q})^{\frac{q}{2-q}}$, from which equality it follows that
\bel{IV-32}
u_0=(V_0C_3^{1-q})^{-\frac{q}{q-2}}V_0^{\frac{2}{q-2}}C_1C_3^{\frac{2}{q-2}}=V_0^{-1}C_1C_3^{-\frac{q(q-1)+2}{q-2}}.
\ee
This implies that $u_0$ can take any value.
\qeda
\mysection{Behaviour of solutions in an exterior domain}
In this section we consider radial solutions of $(\ref{I-1})$ defined in $B^c_{r_0}$. If $u$ is such a solution, it satisfies $(\ref{III-3})$ that we recall 
\bel{V-1}
e^{u(r)}\leq Cr^{-\min\{2,q\}}\qquad\text{for all }\,r\geq 2r_0.
 \ee
Equivalently
\bel{V-2}
e^{u(r)}\leq C\left\{\BA{lll}r^{-2}&\quad\text{if }\,q>2\\
r^{-q}&\quad\text{if }\,1<q<2.
\EA\right.
 \ee
 Applying \rth{sup-2} in $B_{\frac{|x|}{2}}(x)$ for $x\in B_{r_0}^c$ we obtain
 \bth{sup-1*} Let $N\geq 1$ and $q>1$. If $u\in C(B_{r_0}^c)$ is a radial solution of $(\ref{I-1})$ in $B^c_{r_0}$, there holds  in $B^c_{2r_0}$,
  for some positive constant $\tilde C=\tilde C(N,q,M)$,
 \bel{V-3}
|u_r(r)|\leq \tilde C\left\{\BA{lll}r^{-\frac{1}{q-1}}&\quad\text{if }\,q>2\\
r^{-1}&\quad\text{if }\,1<q<2.
\EA\right.
 \ee
 \es
 We will see see later on that in both case the estimate is of the form
 \bel{V-3*}
 |u_r(r)|\leq \tilde C r^{-1}.
 \ee

  \bth{inftyq>2} Let $N\geq 3$ and $q>2$. If $u\in C(B_{r_0}^c)$ is a radial solution of $(\ref{I-1})$ in $B^c_{r_0}$, it satisfies
 \bel{V-4}
\lim_{r\to\infty}r^2e^{u(r)}=2(N-2)\quad\text{and }\;\lim_{r\to\infty}ru_r(r)=-2.
 \ee
 \es
 \Proof As in \rth{posi} we use the systems $(\ref{II-3})$ and $(\ref{II-9})$ with variables variables $x(t)=r^2e^{u(r)}$, $\Phi=-ru_r(r)$  and $\Gth(t)=e^{(2-q)t}$. Here, since $q>2$ and $t=\ln r\to\infty$, we have that $\Gth(t)\to 0$ when $t\to\infty$. The system  in the variable $t$  admits  $O=(0,0,0)$ and 
 $P_0=(2(N-2),2,0)$ for equilibria.  
 The  eigenvalues of the linearised operator at $P_0$ given by $(\ref{Ex-VP})$ are $\gl_1=2-q<0$ and $\gl_2,\gl_3$ are negative too if $N\geq 11$, double with value $2-N$ if $N=10$ or non-real with negative real part $\CR e(\gl_j)=2-N$ if $3\leq N\leq 9$. Therefore $P_0$ is a sink with a domain of attraction which contains some ball $B_a(P_0)$. The eigenvalues of the linearised operator at $O$ 
 $$\left\{\BA{lll}
 x_t=2x\\
 \Phi_t=x-(N-2)\Phi\\
 \Gth_t=(2-q)\Gth,
 \EA\right.$$
 are $2$, $2-q$ and $2-N$, hence $O$ is a saddle point with a $2$-dimensional stable manifold $\CM$ and a unstable trajectory $\CT_0$. \\
 {\it We claim that all the solutions in $B_{r_0}^c$ belong to the domain of attraction of $P_0$ and behave as in $(\ref{V-4})$}.
 Indeed $(x,\Phi)$ satisfies $(\ref{II-3})$ that is
  $$\left\{\BA{lll}
x_t=(2-\Phi)x\\
\Phi_t=x+(2-N)\Phi-Me^{(2-q)t}|\Phi|^q.
\EA\right.$$
Since $e^{u(r)}\to 0$ from $(\ref{V-2})$, $u$ is decreasing, thus $u_r(r)<0$ (the inequality is strict from the equation),  then $\Phi(t)>0$ and $x(t)$ is bounded from $(\ref{V-2})$. Suppose now that 
$\Phi$ is unbounded.\\
- either $\Phi$ is monotone, thus $\Phi(t)\to\infty$. From $(\ref{II-3})$ we have $\Phi_t<\-\frac{N-2}{2}\Phi$ for $t$ large enough, which implies that $t\mapsto e^{\frac{N-2}{2}t}\Phi(t)$ is decreasing. 
Therefore it is bounded which is contradictory.\\
- or there exists a sequence $\{t_n\}$ tending to infinity such that $\Phi(t_n)$ is a local maximum of $\Phi$ and $\Phi(t_n)\to\infty$. Since $\Phi_t(t_n)=0$ we obtain from $(\ref{II-3})$ that 
$x(t_n)=(N-2)\Phi(t_n)+Me^{(2-q)t_n}\Phi(t_n)^q\to\infty$ which contradicts the boundedness of $x(t)$.\\
Therefore $\Phi(t)$ remains bounded. This implies that the system $(\ref{II-3})$ is an exponentially small perturbation of the system associated to Emden-Chandrashekar equation $(\ref{I-E-C})$
  \bel{V-6}\left\{\BA{lll}
x_t=(2-\Phi)x\\
\Phi_t=x+(2-N)\Phi.
\EA\right.\ee
By \rprop {LR} the omega-limit set of any trajectory of $(\ref{II-9})$ is a compact connected subset invariant for $(\ref{V-6})$. This system is well studied (e.g. \cite{Cha}, \cite{BiVe0}). It is known that the solutions of $(\ref{V-6})$ converges to an equilibrium when $t\to\infty$ and this equilibrium is either 
$(2(N-2),2)$ or $(0,0)$. In the first case we obtain $(\ref{V-4})$. In the second case we use system $(\ref{II-9})$ and we analyse the eventual convergence of a solution $(x,\phi,\Gth)$ to 
$O=(0,0,0)$. This would imply that this trajectory belong to the stable manifold $\CM$, therefore 
  \bel{V-7}|\Phi(t)|=O(e^{(2-N)t}).
\ee
If we plug this estimate into $(\ref{II-3})$, we obtain that 
$$x(t)=x(t_0)e^{\int_{t_0}^t(2-\Phi(s))ds}.
$$
Set 
$A(t_0,t)=e^{-\int_{t_0}^t(\Phi(s))ds}$, then there exists $\gth_1>\gth_2>0$ such that $\gth_2\leq A(t_0,t)\leq \gth_1$, for all $t>t_0>0$. This implies 
$$|x(t)|\geq |x(t_0)|e^{2(t-t_0)}\gth_2\quad\text{for all }\,t\geq t_0.
$$
Consequently, if the omega limit set of the trajectory $(x(t),\Phi(t))$ of $(\ref{II-3})$ is $(0,0)$ there must hold $x(t_0)=0$. As a consequence this trajectory must be contained in the plane 
in the plane $x=0$. In such a case $\Phi$ satisfies the equation
  \bel{V-8}\Phi_t=(2-N)\Phi-Me^{(2-q)t}|\Phi|^q.
\ee
This equation is explicitely integrable and we obtain
$$|\Phi(t)|^{-q}\Phi(t)=t^{(q-1)(N-2)}|\Phi(1)|^{-q}\Phi(1)+\frac{M}{(q-1)(N-2)}(t^{(q-1)(N-2)}-t)
$$
{\it This relation is not compatible with $(\ref{V-7})$}. Hence the omega limit set of the trajectory is reduced to $P_0$, which ends the proof.\qeda\medskip

In the case $1<q<2$ we show that all the solutions behave at infinity like the solution of the eikonal equation.
 \bth{inftyq<2} Let $N\geq 3$ and $1<q<2$. If $u\in C(B_{r_0}^c)$ is a radial solution of $(\ref{I-1})$ in $B^c_{r_0}$, it satisfies
 \bel{V-8*}
\lim_{r\to\infty}r^qe^{u(r)}=Mq^q\quad\text{and }\;\lim_{r\to\infty}ru_r(r)=-q.
 \ee
 Note that this behaviour is satisfied both by any regular solution or by the singular solution we have constructed.
 \es
 \Proof We use the system $(\ref{II-4})$ with variables $X(t)=r^qe^{u(r)}$ and $\Phi(t)=-ru_r(r)$, always with $t=\ln r$, that is
$$
\left\{\BA{lll} 
X_t=X(q-\Phi)\\
\Phi_t=(2-N)\Phi-e^{(2-q)t}\left(-M|\Phi|^q+X\right).
\EA
\right.
$$
  Note that $u(r)$ cannot have local minimum, and since  it tends to $-\infty$ by 
 $(\ref{V-2})$ it is decreasing, thus $\Phi(t)\geq 0$. 
Here $X$ is bounded by $(\ref{V-2})$ and also $\Phi$ is bounded by $(\ref{V-3})$. \smallskip

 \nind {\it We claim that $X$ is bounded from below.} We still use the function $G$ defined at $(\ref{IV-G})$, that is 
 $$G(t)=Mq^q\frac{X^2}{2}-\frac{X^3}{3}+e^{(q-2)t}\left((N-2)q\frac{X}{2}-\frac{X_t^2}{2}\right).
 $$
 Then from $(\ref{IV-11*})$, 
 $$G_t(t)=MX^{1-q}\left(qX-Y\right)\left(q^qX^q-Y^q\right)+e^{(q-2)t}\Psi(t),
 $$
 where
 $$\Psi(t)=(q-2)\left((N-2)q\frac{X}{2}-\frac{X_t^2}{2}\right)+\frac{(N-2)q}{2}X_t-\frac{X^3_t}{X}+(N-2)(X_t-qX)X_t.
 $$
 Since $X$ and $\Phi$ are bounded  the same holds with $X_t$ from the equation $(\ref{II-4})$, and because $\frac{X^3_t}{X}=X^2(q-\Phi)^3$, we obtain that $\Psi$ is also bounded by some constant $C>0$.
 Because $\left(qX-Y\right)\left(q^qX^q-Y^q\right)\geq 0$ we obtain that $G_t\geq -Ce^{(q-2)t}$ from what it follows that the function 
 $t\mapsto G(t)+\frac{C}{q-2}e^{(q-2)t}$ is increasing. Since it is bounded, it converges to some limit $\ell$ when $t\to\infty$. Hence the function 
 $t\mapsto Mq^q\frac{X^2(t)}{2}-\frac{X^(t)}{3}$ tends to $\ell$ which implies that $X(t)$ admits a limit $\gl$ when $t\to\infty$,  and $\gl$ satisfies
$Mq^q\frac{\gl^2}{2}-\frac{\gl^3}{3}=\ell$.  If $\Phi(t)$ is not monotone, then at each extremum $t_n$ of $\Phi(t)$ we have 
$$X(t_n)-M\Phi^q(t_n)=(N-2)e^{(q-2)t_n}\Phi(t_n)\Longrightarrow \lim_{n\to\infty}\Phi^q(t_n)=\frac\gl M.
$$ 
This implies that 
$$\limsup_{t\to\infty}\Phi(t)=\left(\frac{\gl}{M}\right)^\frac 1q=\liminf_{t\to\infty}\Phi(t)=\lim_{t\to\infty}\Phi(t).
$$
If $\Phi$ is monotone, then it admits also a limit. In any case we set $\dsps L=\lim_{t\to\infty} \Phi(t)$. From the equation $(\ref{II-4})$ we have
$$\frac{d}{dt}\left(e^{(N-2)t}\Phi(t)\right)=e^{(N-q)t)}(X(t)-M\Phi^q(t)).
$$
Hence 
 \bel{V-9}\Phi(t)=e^{(2-N)(t-t_0)}\Phi(t)+e^{(2-N)t}\int_{t_0}^te^{(N-q)s}(X(s)-M\Phi^q(s))ds.
\ee
If $\gl-ML^q\neq 0$ we obtain a contradiction, since the integration of $(\ref{V-9})$ implies
$$\Phi(t)=\frac{e^{(2-q)t}}{N-q}(\gl-ML^q)(1+o(1)).
$$
Therefore 
$$\lim_{t\to\infty}X(t)=M\left(\lim_{t\to\infty}\Phi(t)\right)^q.
$$
 Again from $(\ref{II-4})$,
 $$X(t)=X(t_0)e^{\int_{t_0}^t(q-\Phi(s)ds}.
 $$
 Then, if $L<q$, $e^{\int_{t_0}^t(q-\Phi(s)ds}\to\infty$, contradiction. If $L>q$, $e^{\int_{t_0}^t(q-\Phi(s)ds}\to 0$, then $X(t)\to 0$, then $\gl=L=0$. We use again the first equation in 
 $(\ref{V-9})$ and we have
 $X_t=\tilde q(t)X(t)$ with $\tilde q=q-\Phi(t)\to q$ when $t\to\infty$. Then for any $t_0>0$ we have
 $$X(t)=X(t_0)e^{\int_{t_0}^t\tilde q(s)ds}.
 $$
 Thus $X(t)\to\infty$, contradiction. Hence 
  \bel{V-10}
  \lim_{t\to\infty}\Phi(t)=q\quad\text{and }\lim_{t\to\infty}X(t)=Mq^q.
\ee
This ends the proof.\qeda

\mysection{Globality of the solutions when $1<q<2$}
\bth{globsub} Let $1<q<2$ and $N\geq 1$. Then the maximal interval of a solution $u$ of $(\ref{II-1})$ is $(0,\infty)$ if $N\geq 2$ and 
$(-\infty,\infty)$ if $N=1$. In that case  $u$ is  symmetric with respect to some $a\in\BBR$.
\es
\Proof For $N\geq 1$ we consider any solution $u$ defined on a maximal open interval $I_u=(\gr,\eta)$ and let $r_0\in I_u$. From Remark 3.2* if $u_r(r_0)\leq 0$, then $\eta=\infty$. If $u_r(r_0)>0$, either $u$ has a unique maximum at some point $r_1\in (r_0,\eta)$, or  $u_r(r_0)> 0$. 
\smallskip

\nind (i) We first prove that $\eta=\infty$. From Remark 3.2*, if $u_r(r_0)\leq 0$, then $\eta=\infty$. If $u_r(r_0)> 0$ has a maximum at $r_1$, then $u$ is decreasing on $(r_1\eta)$ and thus $\eta=\infty$ by the same remark.\\
Hence $u_r>0$ on $(r_1,\eta)$) and we can encounter two possibilities: either $\eta=\infty$ and  $u(r)\to L\in (0,\infty)$ when 
$r\to\infty$, which contradicts the upper estimate $(\ref{III-3})$, or $\eta<\infty$, and we are left with this case, with again two possibilities:\\
If $\dsps \lim_{r\to\eta}u(r)=L<\infty$, then 
$$(-r^{N-1}u_r)_r=r^{N-1}(e^u-u_r^q)\leq k=\eta^{N-1}e^L.
$$
Therefore the function $r\mapsto r^{N-1}u_r+kr$ is increasing and it tends to $\infty$ since the solution is the maximal one and $u$ remains bounded. This implies that $u_r(r)\to\infty$ and $e^{u(r)}=o(u_r^q)$ when $r\to\eta$. Hence the equation becomes $-u_{rr}+Mu_r^q(1+o(1))$. Setting $v=u_r$ we obtain by integration that 
$$v(r)=u_r(r)=\left(M(q-1)(\eta-r)\right)^{-\frac{1}{q-1}}(1+o(1))\quad\text{as }r\to\eta,
$$
but this is not possible since $1<q<2$ and thus $u_r$ is not integrable near $\eta$, which contradicts $L<\infty$. 
Consequently there must hold $L=\infty$. Again $u_r$ is monotone near $\eta$ since at any extremal point $r$ where 
$u_{rr}(r)=0$, we have that $u_{rrr}(r)=\left(\frac{N-1}{r^2}-e^{u}\right)u_r$. Since $\frac{N-1}{r^2}-e^{u}\to-\infty$ when $r\to\eta$, then $u_{rrr}$ keeps a constant sign. Therefore $u_r$ is non-decreasing and clearly
$$u_{rr}=Mu_r^q-\frac{N-1}{r}u_r-e^{u}\geq 0.
$$
Then $e^{u}\leq Mu_r^q(1+o(1))$, and by integration we deduce that $e^{u}\leq C(\eta-r)^{-q}$. Now we set $u(r)=v(h(r))$ where $h(r)=\eta-r$ and consider the system in the variables
\bel{G-1}
X=h^qe^{v},\quad \Phi=-hv_r(h),\quad \,h=e^\gt,
\ee
which is
\bel{G-2}\left\{\BA{lll}
X_\gt=X(q-\Phi)\\
\Phi_\gt=\Phi
+e^{(2-q)\gt}(X-M\Phi^q)-\frac{N-1}{\eta-e^\gt}\Phi,
\EA\right.
\ee
and $\gt\to-\infty$ when $r\to\eta$. We claim now that $\Phi$ is bounded. If it is not the case and $\Phi(\gt)$ is not monotone, at any $\tilde\gt$ where $\Phi_\gt(\tilde\gt)=0$ we obtain
$$\Phi(\tilde\gt)(1+o(1))=e^{(2-q)\tilde\gt}(M\Phi^q(\tilde\gt)-X(\tilde\gt)\leq e^{(2-q)\tilde\gt}M\Phi^q(\tilde\gt),
$$
and thus $M\Phi^{q-1}(\tilde\gt)\geq e^{-(2-q)\tilde\gt}(1+o(1))$. Since this is valid in particular at the local minima 
of $\Phi$, we deduce that $\Phi(\gt)\to\infty$ when $\gt\to-\infty$, with the additional information that  $\Phi_{\gt\gt}(\tilde\gt)\geq 0$. Since 
$$ \Phi_{\gt\gt}(\tilde\gt)=e^{(2-q)\tilde \gt}\left(2X(\tilde\gt)-(X\Phi)(\tilde\gt)-M(2-q)\Phi^q(\tilde\gt)\right)
-\frac{N-1}{\eta-e^{\tilde\gt}}e^{\tilde\gt}\Phi(\tilde\gt)\left(1+\frac{e^\gt}{\eta-e^{\tilde\gt}}\right)\geq 0,
$$
which yields
$$2X(\tilde\gt)-(X\Phi)(\tilde\gt)-M(2-q)\Phi^q(\tilde\gt)\geq e^{(q-1)\tilde\gt}\Phi(\tilde\gt)\left(1+\frac{e^\gt}{\eta-e^{\tilde\gt}}\right), 
$$
which in turn implies that $\Phi(\tilde\gt)$ is bounded, contradiction. As a consequence $\Phi$ is monotone and tends to $\infty$ at $-\infty$. From the first equation in $(\ref{G-2})$, $X=-X\Phi(1+o(1))$ when $\gt\to-\infty$, which is impossible if $X$ remains bounded. As a consequence the function $\Phi$ is bounded and the system $(\ref{G-2})$ is an exponential perturbation of the system 
\bel{G-3}\left\{\BA{lll}
X_\gt=X(q-\Phi)\\
\Phi_\gt=\Phi.
\EA\right.
\ee
This implies that $(X(\gt),\Phi(\gt))$ converges to the unique stationary point of $(\ref{G-3})$ which is $(0,0)$.
\\
Returning to the variable $u$, we get that $u_r=o((\eta-r)^{-1})$ near $r=\eta$, then 
$$u_{rr}=M|u_r|^q-\frac{N-1}{r}u_r-e^{u}=o((\eta-r)^{-q}),
$$
which implies that $u_r(r)=o((\eta-r)^{1-q})$. Thus $u_r$ is integrable near $\eta$ and $u(r)$ admits a finite limit when
$r\to\eta$, which is a contradiction. As a consequence we deduce that:\smallskip

\nind (ii) If $N=1$ then the solution is defined on whole $\BBR$, indeed $r\mapsto u(-r)$ is also a solution, thus $I_u=(\gr,\eta)=(-\infty,\infty)$.\smallskip

\nind (iii) If $N\geq 2$, we claim that $\gr=0$. We proceed again by contradiction, assuming that $\gr>0$. Clearly 
$u$ is monotone near $\gr$.\\
- If $u$ is increasing near $\gr$ it has a finite limit $L$ at $\gr$ because 
\bel{G-4}-(r^{N-1}u_r)_r=r^{N-1}(e^u-Mu_r^q)\leq C\quad\text{as }r\to\gr,
\ee
which implies that $r^{N-1}u_r+Cr$ is increasing. Then $u_r$ admits a finite limit at $\gr$. Combined with the fact that 
$u(r)\to L$, we see that $\gr$ cannot be the infimum of $I_u$. Consequently $u(r)$ decreases to $-\infty$ when $r\to\gr$. In such a case $e^{u(r)}\to 0$. Moreover $u_r$ is also monotone near $\gr$: indeed at each local extremum  
$\tilde r$ of $u_r$ we have $u_{rr}(\tilde r)=0$, we have 
$$u_{rrr}(\tilde r)=\left(\frac{N-1}{\tilde r^2}\right)u_r(\tilde r)\to \frac{N-1}{\gr^2}>0 \quad\text{as }\tilde r\to\gr.
$$
As a consequence these local extrema are local minima of $u_r$; necessary $u_r(r)$ cannot oscillate and it tends to 
$\infty$ when $r\to\gr$. Again this fact implies that $u_{rr}=Mu_r^q(1+o(1))$ and by integration we encounter a contradiction since $1<q<2$.\\
- If $u$ is decreasing near $\gr$, then $u(r)\to L<\infty$ when $r\to\gr$. By $(\ref{G-4})$, $-(r^{N-1}u_r)_r\leq C$, hence 
$r\mapsto -r^{N-1}u_r+Cr$ is increasing and thus $\dsps\lim_{r\to\gr}u_r(r)=-\infty$ since $\gr$ is an end point $I_u$. We have again that $u_r$ is monotone and $u$ is convex. Therefore 
$0\leq u_{rr}=M|u_r|^q-\frac{N-1}{r}u_r-e^{u}$. hence $e^u\leq M|u_r|^q(1+o(1))$ and by integration we deduce that
$e^{u(r)}\leq C(r-\gr)^{-q}$. As in (i) we set $h(r)=r-\gr$  $e^\gt=r-\gr$ and deine $X(\gt)$ and $\Phi(\gt)$ as in $(\ref{G-1})$. Since $(X,\Gf)$ satisfies 
$$\left\{\BA{lll}
X_\gt=X(q-\Phi)\\
\Phi_\gt=\Phi
+e^{(2-q)\gt}(X-M\Phi^q)-\frac{N-1}{\gr+e^\gt}\Phi,
\EA\right.
$$
we obtain a contradiction as in this case. This ends the proof.\qeda
\mysection{Appendix}
In \rlemma{dev1} formula $(\ref{IV-13})$ gives an  expansion of $u$ near zero at the order $1$. We show below that $u(r)$ satisfies an {\it unique} expansion of order $n$ of the form
   \bel{IV-20}\BA{lll}\dsps
u(r)=\ln \frac {Mq^q}{r^q}+a_1r^{q-2}+a_2r^{2(q-2)}+...+r^{n(q-2)}\left(a_n+o(1)\right)
\text{as }r\to 0,
\EA\ee
where the $a_n$ can be computed by induction.

\bth . Let $N\geq 2$ and $q>2$. There exists a unique sequence of real numbers
$\{a_n\}_{n\geq 1}$ such that any radial function $u$ satisfying $(\ref{VI-1})$ in $B_{r_0}\setminus\{0\}$ and  $(\ref{IV-4})$ there holds   for any $n\in\BBN^*$,
  \bel{VI-2}
  u(r)=\ln\frac{q^q}{r^q}+a_1r^{q-2}+a_2r^{2(q-2)}+...+a_nr^{n(q-2)}(1+d_n(r)),
\ee
where $d_n(r)\to 0$ when $r\to 0$.
\es
\Proof For the sake of simplicity we consider equation $(\ref{II-1})$ with $M=1$. Actually, this is not a restriction since if $u^*$ satisfies $(\ref{I-1})$, the function 
$x\mapsto u(x)=u^*(M^{\frac{1}{q-2}}x)+M^{\frac{2}{q-2}}$ satisfies
  \bel{VI-1}-\Gd u+|\nabla u|^q-e^{u}=0.
\ee
Since any solution is decreasing, we consider $u_r$ as a function of $u$ and set
$$u_r=-f(u).
$$
Then $u_{rr}=-f(u)\frac{df}{du}=f^q(u)-e^{u}+\frac{N-1}{r}f(u)$. Thus we obtain the system
  \bel{VI-3}\left\{\BA{lll}\dsps
\frac{df}{du}=\frac{N-1}{r}+\frac{f^q(u)-e^{u}}{f(u)}\\[3mm]
\dsps \frac{dr}{du}=-\frac{1}{f(u)}.
\EA\right.\ee
Using $(\ref{IV-4})$ we have
$$f(u)=e^{\frac uq}(1+o(1))\,\text{ and }\;r(u)=qe^{-\frac uq}(1+o(1))\,\text{ as }\;u\to \infty.
$$
We set
  \bel{VI-4}
\varpi=e^{-\frac uq}f(u),\quad V=\varpi-1\;\text{ and }\;\gth=\frac{q-2}{q}.
\ee
Notice that $e^{\gth u}\to\infty$ since $q>2$. The new system in $(\varpi,r)$ is
  \bel{VI-5}\left\{\BA{lll}\dsps
\frac{d\varpi}{du}=-\frac{\varpi}{q}+\frac{N-1}{re^{\frac uq}}+\frac{\varpi^q-1}{\varpi}e^{\gth u}\\[3mm]
\dsps \frac{dr}{du}=-\frac{e^{-\frac uq}}{\varpi}.
\EA\right.\ee
{\it Step 1: Development of order 1}. We look for $A_1$ such that 
$$\varpi(u)=1+A_1e^{-\gth u}(1+o(1))\,\text{ as }\;u\to \infty.
$$
Equivalently $V(u)=A_1e^{-\gth u}(1+o(1))$. The system in $V$ and $r$ is 
  \bel{VI-6}\left\{\BA{lll}\dsps
\frac{dV}{du}=-\frac{1+V}{q}+\frac{N-1}{re^{\frac uq}}+\frac{(1+V)^{q}-1}{1+V}e^{\gth u}\\[3mm]
\dsps \frac{dr}{du}=-\frac{e^{-\frac uq}}{1+V}.
\EA\right.\ee
We set $\psi(u)=e^{\gth u}V(u)$, then $\psi(u)=o(e^{\gth u})$ when $u\to\infty$. Hence
$$e^{-\gth u}\left(\frac{d\psi}{du}-\frac{q-2}{q}\psi\right)=-\myfrac{1+e^{-\gth u}\psi}{q}+\frac{N-1}{re^{\frac uq}}+e^{\gth u}\frac{(1+e^{-\gth u} )^q\psi-1}{1+e^{-\gth u}\psi},
$$
$$\frac {dr}{du}=\frac{e^{-\frac uq}}{1+e^{-\gth u}\psi}.
$$
As a consequence
$$\BA{lll}\dsps\frac{d\psi}{du}=\frac{q-2}{q}\psi-\frac{e^{\gth u}}{q}-\frac \psi q+\frac{N-1}{re^{\frac uq}}e^{\gth u}+qe^{\gth u}\psi(1+o(1))\\[3mm]
\phantom{\dsps\frac{d\psi}{du}}\dsps=qe^{\gth u}\psi(1+o(1)))+e^{\gth u}\left(-\frac 1q+\frac{N-1}{q}(1+o(1))\right),
\EA$$
that we can write under the following form
  \bel{VI-7}\frac{d\psi}{du}-qe^{\gth u}\psi(1+o(1))=\frac{N-2}{q}e^{\gth u}(1+o(1)).
\ee
- Either $\psi$ is not monotone for large $u$. At each $u$ where $\frac{d\psi}{du}=0$, we get that $\psi(u)=\frac{2-N}{q^2}(1+o(1))$. These implies 
$\dsps\lim_{u\to\infty}\psi(u)=\frac{2-N}{q^2}$. \\
- Or $\psi$ is monotone and there holds $\dsps\lim_{u\to\infty}\psi(u)=\pm\infty$. Therefore $\frac{d\psi}{du}=qe^{\gth u}\psi(1+o(1))$, which implies $\frac{d}{du}(\ln |\psi|-e^{(q-2)u})>0$ and thus 
$|\psi(u)|\geq Ke^{e^{(q-2)u}}$. Since $\psi(u)=o(e^{\gth u})$ we obtain a contradiction.\smallskip

\nind- Or $\psi$ is monotone and $\psi(u)\to \ell$ when $u\to\infty$ for some $\ell$. If $\ell\neq \frac{2-N}{q^2}$, we obtain the following relation  
$\frac{d\psi}{du}=\left(q\ell-\frac{N-2}{q}\right)(1+o(1)))$ which is not compatible with $\psi(u)\to \ell$. \\
In all the cases we obtain
  \bel{VI-8}\lim_{u\to\infty}\psi(u)=\frac{2-N}{q^2}.
\ee
This yields the first term of our expansion
  \bel{VI-9}
  \varpi=1+V=1+A_1e^{-\gth u}(1+o(1))\quad\text{with }\;A_1=\frac{2-N}{q^2}.
\ee
as a consequence
  \bel{VI-10}
\frac{dr}{du}=-\frac{e^{-\frac uq}}{1+A_1e^{-\gth u}(1+o(1))}=-e^{-\frac uq}\left(1-A_1e^{-\gth u}(1+o(1))\right)=-e^{-\frac uq}+A_1e^{\frac {1-q}q u}(1+o(1)).
\ee
Integrating this relation and using that $r(u)\to 0$ when $u\to\infty$, we obtain
  \bel{VI-11}\BA{lll}\dsps
r(u)=-\int_u^\infty\frac{dr}{ds}ds=\int_u^\infty e^{-\frac sq}ds-A_1\int_u^\infty e^{\frac {1-q}{q}s}(1+\ge(s))ds\\[3mm]
\phantom{\dsps
r(u)}\dsps=qe^{-\frac uq}+\frac{qA_1}{1-q}e^{\frac{1-q}{q}u}(1+o(1))=qe^{-\frac uq}+\frac{N-2}{q(q-1)}e^{\frac{1-q}{q}u}(1+o(1)),
\EA
\ee
since for any $\ge_0>0$and $|\ge(s)|\leq\ge_0$ there exists $u_{\ge_0}$ such that for any $u\geq u_{\ge_0}$ there holds
$$\left|\int_u^{\infty}e^{\frac{1-q}qs}(1+\ge(s))ds-\int_u^{\infty}e^{\frac{1-q}qs}ds\right|\leq \frac{q\ge_0}{q-1}.
$$
{\it Step 2: Development of order n}. We proceed by induction assuming that we have already obtained the development of $\varpi$ at the order $n-1$
  \bel{VI-12}
\varpi(u)=1+V=1+A_1e^{-\gth u}+A_2e^{-2\gth u}+...+e^{-(n-1)\gth u}\left(A_{n-1}+\ge_{n-1}(u)\right),
\ee
where the $A_j$ depend on $N$ and $q$ and $\ge_{n-1}(u)\to 0$ when $u\to\infty$. Consequently 
\bel{VI-bis}\BA{lll}\dsps
\frac{dr}{du}=-\frac{e^{-\frac uq}}{1+V}=-\frac{e^{-\frac uq}}{1+A_1e^{-\gth u}+A_2e^{-2\gth u}+...+e^{-(n-1)\gth u}\left(A_{n-1}+\ge_{n-1}(u)\right)}\\
[3mm]\phantom{\dsps
\frac{dr}{du}}\dsps
=e^{-\frac uq}\left(1+B_1e^{-\gth u}+B_2e^{-2\gth u}+...+e^{-(n-1)\gth u}\left(B_{n-1}+\tilde \ge_{n-1}(u)\right)\right),
\EA
\ee
where $\tilde \ge_{n-1}(u)\to 0$ as $u\to\infty$, with $B_1=-A_1$   and the other coefficients can be made explicit through a lengthy but explicit computation. Set
\bel{VI-ter}\varpi=1+A_1e^{-\gth u}+A_2e^{-2\gth u}+...+A_{n-1}e^{-(n-1)\gth u}+e^{-n\gth u}\Phi.
\ee
We do not know if $\Gf$ is bounded but only that $\Phi=\tilde\ge_{n-1}e^{\gth u}$, thus $e^{-\gth u}\Phi(u)\to 0$ when $u\to\infty$. Thus
\bel{help0}\frac{dV}{du}=-\gth A_1e^{-\gth u}-2\gth A_2e^{-2\gth u}-...-(n-1)\gth A_{n-1}e^{-(n-1)\gth u}-n\gth \Phi e^{-n\gth u}+\frac{d\Phi}{du}e^{-n\gth u}.
\ee
We recall the expression
\bel{help1}\frac{dV}{du}=-\frac{1+V}{q}+\frac{N-1}{re^{\frac uq}}+e^{\gth u}\frac{(1+V)^q-1}{1+V},
\ee
obtained by using the following expansion at the order $n-1$ (actually, valid at any order)
\bel{help1-bis}\frac{(1+V)^q-1}{1+V}
=qV\left(1+g_1V+...+g_kV^k+...+g_{n-1}V^{n-1}(1+\eta_{n-1}(V))\right),
\ee
where $g_1=\frac{q-3}{2}$ and the $g_k$ are polynomials in $k$ and $q$ and also $\eta_{n-1}(V)\to 0$ when $V\to 0$.\\ Moreover, by integration of $(\ref{VI-bis})$, we obtain
$$\BA{lll}\dsps
r(u)=-\int_u^\infty\frac{dr}{ds}ds=-\int_u^\infty \!\!\!e^{-\frac sq}ds+B_1\int_u^\infty \!\!\!e^{\frac{1-q}{q}s}+...+\int_u^\infty \!\!\!e^{-\frac{1+(n-1)(q-2)}{q}s}(B_{n-1}+\hat \ge_{n-1}s))ds\\[4mm]
\phantom{r(u)}\dsps=qe^{-\frac uq}\left(1+\frac{B_1}{q-1}e^{-\gth u}+...+\frac{B_{n-1}}{1+(n-1)(q-2)}e^{-(n-1)\gth u}+\hat\ge_{n-1}(u)e^{-(n-1)\gth u}\right),
\EA$$
therefore
$$\BA{lll}\dsps
\frac{N-1}{re^{\frac uq}}=\frac{N-1}{q}\frac{1}{1+\frac{B_1}{q-1}e^{-\gth u}+...+\frac{B_{n-1}}{1+(n-1)(q-2)}e^{-(n-1)\gth u}+\tilde\ge_{n-1}(u)e^{-(n-1)\gth u}}\\[3mm]
\phantom{\dsps
\frac{N-1}{re^{\frac uq}}}
\dsps=\frac{N-1}{q}\left(1+C_1e^{-\gth u}+...+C_{n-1}e^{-(n-1)\gth u}+\ge^*_{n-1}(u)e^{-(n-1)\gth u}\right),
\EA$$
where $C_1=-\frac B{q-1}$ and $C_k$ is a polynomial in $q$ and $k$. Then using $(\ref{help1})$ and $(\ref{help1-bis})$,
\bel{help2}\BA{lll}\dsps
\frac{dV}{du}=-\gth A_1e^{-\gth u}-2\gth A_2e^{-2\gth u}...-(n-1)\gth A_{n-1}e^{-(n-1)\gth u}-n\gth e^{-n\gth u}\Phi+e^{-n\gth}\frac{d\Phi}{du}\\[3mm]
\phantom{\dsps
\frac{dV}{du}}\dsps =\frac{N-1}{q}\left(1+C_1e^{-\gth u}+...+C_{n-1}e^{-(n-1)\gth u}+\ge^*_{n-1}(u)e^{-(n-1)\gth u}\right)\\[3mm]
\phantom{-------}\dsps -\frac{1+V}{q}+qVe^{\gth u}\left(1+g_1V+...+g_kV^k+...+g_{n-1}V^{n-1}(1+\eta_{n-1}(V))\right),
\EA\ee
and, by the definition of $\Phi$, from $(\ref{VI-ter})$ there holds 
\bel{help3}\BA{lll}Ve^{\gth u}=A_1+A_2e^{-\gth u}+...+A_{n-1}e^{-(n-2)\gth u}+e^{-(n-1)\gth u}\Phi, 
\EA\ee
and 
\bel{help4}
V=A_1e^{-\gth u}+A_2e^{-2\gth u}+...+A_{n-1}e^{-(n-1)\gth u}+e^{-(n-1)\gth u}\tilde\ge_{n-1}(u).
\ee
Then, using $(\ref{help3})$ we compute the expression $1+g_1V+g_2V^2+...+g_nV^{n}(1+\eta_{n-1}(V))$.
Since $|\ge_{n-1}(u)|<\!<1$, we write for any $k=1,...,n-1$,
$$V^k=(a+\gt)^k\,\text{ with }\,a=A_1e^{-\gth u}+A_2e^{-2\gth u}+...+A_{n-1}e^{-(n-1)\gth u}\,\text{ and }\,\gt =e^{-(n-1)\gth u}\tilde \ge_{n-1}(u).
$$
Then
$$\BA{lll}
|(a+\gt)^k-a^k|\leq k|\gt|(|a|+|\gt|)^k\leq k|\gt|\left(|A_1|+...+|A_{n-1}|+|\tilde\ge_{n-1}(u)|\right)^{k-1}|\tilde\ge_{n-1}(u)|\\[2mm]
\phantom{|(a+\gt)^k-a^k|}
\leq k|A_{n-1}|\left(|A_1|+...+|A_{n-1}|+1\right)^k
e^{-(n-1)u}|\tilde\ge_{n-1}(u)|\\[2mm]
\phantom{|(a+\gt)^k-a^k|}=c_k|\tilde\ge_{n-1}(u)|e^{-(n-1)u};
\EA$$
now
$$g_kV^k=g_k\left(A_1e^{-\gth u}+A_2e^{-2\gth u}+...+A_{n-1}e^{-(n-1)\gth u}\right)^k+\gd_{n-1,k}(u),$$
 with $|\gd_{n-1,k}(u)|\leq c_k|g_k||\ge_{n-1}(u)|$. Therefore
$$\BA{lll}
\dsps 1+g_1V+...+g_kV^k+...+g_{n-1}V^{n-1}(1+\eta_{n-1}(V))\\[3mm]
\dsps=1+g_1\left(A_1e^{-\gth u}+A_2e^{-2\gth u}+...+A_{n-1}e^{-(n-1)\gth u}\right)+...\\[3mm]
\phantom{----------}+g_k\left(A_1e^{-\gth u}+A_2e^{-2\gth u}+...+A_{n-1}e^{-(n-1)\gth u}\right)^k\\[3mm]
\phantom{----------}
+g_n\left(A_1e^{-\gth u}+A_2e^{-2\gth u}+...+A_{n-1}e^{-(n-1)\gth u}\right)^n+e^{-(n-1)\gth u}\gd_{n-1}(u),
\EA$$
then 
\bel{VI-n}\BA{lll}
1+g_1V+...+g_kV^k+...+g_{n-1}V^{n-1}(1+\eta_{n-1}(V))\\[2mm]
\phantom{-----}=1+D_1e^{-\gth u}+D_2e^{-2\gth u}+...+D_{n-1}e^{-(n-1)\gth u}+e^{-(n-1)\gth u}\gd_{n-1}(u)
\EA\ee
with $|\gd_{n-1}(u)|\leq C_{n-1}|\tilde\ge_{n-1}(u)|$ and the coefficients $D_j$ can be explicitely computed. Next we compute from the expression of $\varpi$, $(\ref{help0})$-$(\ref{help4})$ and $(\ref{VI-n})$,
\bel{help4bis}\BA{lll}
\dsps 
e^{-n\gth}\left(\left(\frac 1q-n\gth\right)\Phi+\frac{d\Phi}{du}\right)\\[4mm]
\phantom{----}\dsps =\frac{dV}{du}+\gth A_1e^{-\gth u}+2\gth A_2e^{-2\gth u}+...+(n-1)\gth A_{n-1}e^{-(n-1)\gth u}+\frac 1qe^{-n\gth}\Phi\\[3mm]
\phantom{----}\dsps
=\frac{dV}{du}+\left(\gth-\frac 1q\right)A_1e^{-\gth u}+\left(2\gth-\frac 1q\right)A_2e^{-2\gth u}+...+\left((n-1)\gth-\frac 1q\right)A_{n-1}e^{-(n-1)\gth u}\\[3mm]
\phantom{----}\dsps
=\frac{N-1}{q}\left(1+C_1e^{-\gth u}+...C_ke^{-k\gth u}+...+\left(C_{n-1}+\tilde\ge_{n-1}(u)\right)e^{-(n-1)\gth u}\right)\\[3mm]
\phantom{-----}\dsps-\frac{1+V}{q}+qVe^{\gth u}\left(1+g_1V+...+g_kV^k+...+g_{n-1}V^{n-1}(1+\eta_{n-1}(V))\right)\\[3mm]
\phantom{----}\dsps
=\frac{N-1}{q}\left(1+C_1e^{-\gth u}+...C_ke^{-k\gth u}+...+\left(C_{n-1}+\tilde\ge_{n-1}(u)\right)e^{-(n-1)\gth u}\right)\\[3mm]
\phantom{-----}\dsps
+q\left(A_1+A_2e^{-\gth u}+...+A_{n-1}e^{-(n-2)\gth u}+e^{-(n-1)}\Phi\right)\\[3mm]
\phantom{------}\dsps
\ti \left(1+D_1e^{-\gth u}+D_2e^{-2\gth u}+...+D_{n-1}e^{-(n-1)\gth u}+e^{-(n-1)\gth u}\gd_{n-1}(u)\right)\\[3mm]
\phantom{-----}
+\left(\gth-\frac 1q\right)A_1e^{-\gth u}+\left(2\gth-\frac 1q\right)A_2e^{-2\gth u}+...+\left((n-1)\gth-\frac 1q\right)A_{n-1}e^{-(n-1)\gth u}\\[3mm]
\phantom{----}\dsps
=\frac{N-2}{q}+qA_1+E_1e^{-\gth u}+E_2e^{-2\gth u}+...\\[3mm]
\phantom{-----}+E_{n-1}e^{-(n-1)\gth u}(1+\gd_{n-1}(u))+qe^{-(n-1)\gth u}\Phi (1+\tilde \gd_n(u)),
\EA\ee
where $\tilde \gd_n(u)\to 0$.\\  
Because of $(\ref{VI-9})$ $\frac{N-2}{q}+qA_1=0$. {\it We claim that $E_2=E_3=...=E_{n-2}=0$}. If this does not hold, let $k\in [1,n-2]$ be the smaller integer such that $E_k\neq 0$.
Then 
$$\left(\frac 1q-n\gth-qe^{\gth u}(1+\tilde\gd_{n}(u))\right)\Phi+\frac{d\Phi}{du}=E_ke^{(n-k)\gth u}.
$$
Since $\Phi(u)=o(e^{\gth u})$, then $\left(\frac 1q-n\gth-qe^{\gth u}(1+\tilde\gd_{n}(u))\right)\Phi=o(e^{2\gth u})$. Thus 
$$\frac{d\Phi}{du}=E_ke^{(n-k)\gth u}(1+o(1)),
$$
which implies $|\Phi|\geq Ce^{(n-k)\gth u}$, contradiction. Therefore we are led to the following relation
\bel{help5}
\left(\frac 1q-n\gth-qe^{\gth u}(1+\tilde\gd_{n}(u)\right)\Phi+\frac{d\Phi}{du}=E_{n-1}e^{\gth u}(1+\tilde \gd_n(u)),
\ee
which yields
\bel{help6}-qe^{\gth u}(1+\gd'_{n}(u))\Phi+\frac{d\Phi}{du}=E_{n-1}e^{\gth u}(1+\tilde \gd_n(u)).
\ee
Finally we conclude as in the Step 1, considering the three possibilities for $\Phi$:\\
- either non-monotone near infinity which implies that $\lim_{u±to\infty}=-\frac{E_{n-1}}{q}:=A_n$,\\
- or 
$\Phi$ is monotone and tends to $\pm\infty$ and this would lead to a contradiction with $\Phi (u)=o(e^{\gth u})$ when $u\to\infty$, or  $\Phi$ is monotone and bounded, thus it admits a limit which is necessarily $-\frac{E_{n-1}}{q}$. Thus we have obtained the development at the order $n$ of the function $\varpi(u)$. This yields a development of 
$\frac{dr}{du}$ by the second equation in system  $(\ref{VI-6})$ under the form 
$$\frac{dr}{du}=e^{-\frac uq}\left(1+B_1e^{-\gth u}+B_2e^{-2\gth u}+...+e^{-n\gth u}\left(B_{n}+\tilde \ge_{n}(u)\right)\right).
$$
From this we obtain the expansion of $r(u)$ by integration and finally obtain algebraically the expansions in the variable $r$
$$e^{\frac{u(r)}{q}}=\frac{q}{r}\left(1+b_1r^{q-2}+...+(b_n+o(1))r^{n(q-2)}\right),
$$
and 
$$u(r)=\ln\frac{q^q}{r^q}+a_1r^{q-2}+...+r^{n(q-2)}(a_n+o(1)).
$$
It is noticeable that all the terms of order $n$ in these expansions are polynomial in $N,n,q$.
\qeda\medskip

\nind\Remark The fact that a solution of eikonal type can be expressed formally by a series in the variable $r^{q-2}$ is a strong presumption  
for uniqueness. However it appears very difficult to prove that the radius of convergence of this series is positive. This type of question is similar 
to the existence and uniqueness problem encountered in the study of the singular solutions of the capillary equation as shown in \cite{CoFi1} and \cite{CoFi2}. However in this problem, the formal series obtained in the study of singular solutions can be proved to have a null radius of convergence. \medskip

\nind\Remark If $N=2$ there exists an explicit solution of $\ref{VI-1}$, namely $u^*(x)=\ln\frac{q^q}{|x|^q}$ and clearly it satisfies $(\ref{IV-4})$. As a consequence all the terms $D_j$ in the expansion of another radial solution of $\ref{VI-1}$ satisfying $(\ref{IV-4})$ are zero, which means that for any $n$, there holds
$$u(r)=\ln\frac{q^q}{r^q}+o(r^{n(q-2)})\quad\text{as }r\to 0.
$$
This situation is reminiscent of a result of Brezis and Nirenberg \cite[Theorem 2]{BrNi} dealing with the equation
$$-\Gd u+|\nabla u|^2=h^2(u).
$$

\begin{thebibliography}{99}

\bibitem{AnLei} Anderson L.R., Leighton W., \textit {Liapunov functions for autonomous systems of second order}.  J. Math. Anal. Appl. {\bf 23} (1968), 645-664.

\bibitem{BiGaVe2} Bidaut-V\'{e}ron M.F., Garcia-Huidobro M., V\'{e}ron L.,
\textit{ A priori estimates for elliptic equations with reaction terms involving the function and its gradient}. Math. Ann. {\bf 378} (2020), 13-56.

\bibitem{BiGaVe3} Bidaut-V\'{e}ron M.F., Garcia-Huidobro M., V\'{e}ron L,
\textit{Radial solutions of scaling invariant nonlinear elliptic equations
 with mixed reaction terms}. Discrete Contin. Dyn. Syst. {\bf 40} (2020), 933-982.
 
 \bibitem{BiGaVe4} Bidaut-V\'{e}ron M.F., Garcia-Huidobro M., V\'{e}ron L,
\textit{Singular solutions of some elliptic equations involving mixed absorption-reaction}. Discrete Contin. Dyn. Syst. {\bf 42} (2022), 3861-3930.

  
 \bibitem{BiVe0} Bidaut-V\'{e}ron M.F., V\'{e}ron L.,\textit{ Nonlinear elliptic equations on compact Riemannian manifolds and asymptotics of Emden equations}, Invent. Math. {\bf106} (1991), 489-539.


\bibitem{BiVe1} Bidaut-V\'{e}ron M.F., V\'{e}ron L.,\textit{ Local behaviour of the solutions of the Chipot–Weissler
equation}, Calc. Var. Partial Differential Equations {\bf 62} (2023), paper No. 241, 59 pp.


 \bibitem{BiVe2} Bidaut-V\'{e}ron, V\'{e}ron L,
\textit{Singularities and asymptotics of solutions of the  Emden-Chandrasekhar-Riccati equation}, arXiv: 2506.01136  [math.AP], submitted.


\bibitem{BrLi} Brezis H., Lions P. L.,\textit{ A note on isolated singularities for linear elliptic equations} in  Mathematical analysis and applications, Part A. 263-266, {\it Adv. Math. Suppl. Stud.} {\bf 7a} (1981).

\bibitem{BrMe} Brezis H., Merle F.,\textit{ Uniform estimates and blow-up behavior for solutions of $-\Gd u=V(x)e^u$ in two dimensions.}
Comm. Partial Differential Equations {\bf 16} (1991), 1223-1253.

\bibitem{BrNi} Brezis H., Nirenberg L.,\textit{ Removable singularities for nonlinear elliptic equations}. Topol. Methods Nonlinear Anal. {\bf 9}
(1997), 201-219.

\bibitem{Ca-Cr} L. Caffarelli, M. Crandall.,\textit{ Distance functions and almost global solutions of eikonal equations}.  Comm. Partial Differential Equations {\bf 35} (2010), 391-414. 

\bibitem{Cha} Chandrasekhar S., Introduction to Stellar Structure. {\it University of Chicago Press}  (1939).

\bibitem{ChWe1} Chipot M., Weissler F., \textit{On the elliptic problem $\Gd u - |\nabla u|^q+ \gl u^p$ = 0}, 
Nonlinear diffusion equations and their equilibrium states, I (Berkeley, CA, 1986), 237-243, Math. Sci. Res. Inst. Publ., {\bf 12}, Springer, New York, (1988).

\bibitem{ChWe2} Chipot M., Weissler F., \textit{Some blowup results for a nonlinear parabolic equation with a gradient term}, SIAM J. Math. Anal. {\bf 20} (1989), 886-907.

\bibitem{CoFi1} Concus P., Finn R., \textit{ A singular solution of the capillary equation. II. Existence.} Invent. Math. {\bf 29} (1975), 143-148.

\bibitem{CoFi2} Concus P., Finn R.,  \textit{ A singular solution of the capillary equation. II. Uniqueness.} Invent. Math. {\bf 29} (1975), 149-159.


\bibitem {Emd} Emden R.  {\it Gaskugeln: Anwendungen der mechanischen W\"armetheorie auf kosmologische und meteorologische Probleme. B}. Teubner, Leipzig (1907).

\bibitem{La-Li} Lasry J. M., Lions P-L., \textit{Nonlinear elliptic equations with singular boundary conditions and stochastic control with
state constraints. I. The model problem}, Math. Ann. {\bf 283} (1980), 583-630.

\bibitem{LiVe} Li Y., V\'{e}ron L., \textit{ Isolated singularities of solutions of a 2-D diffusion equation with mixed reaction}, arXiv:2408.14246v2 [math.AP] (2024), submitted.

\bibitem{Li1} Lions P-L., \textit{Quelques remarques sur les probl\`emes quasilin\'eaires elliptiques du second ordre}, J. Anal. Math. {\bf 45} (1985), 234-254.

\bibitem{LoRy} Logemann H., Ryan E. P., \textit{Non-autonomous systems: Asymptotic behavior and weak
invariance principles}. J.  Diff. Equations {\bf 189} (2003), 440-460.


\bibitem{Ng1} Nguyen P. T., \textit{Isolated singularities of positive solutions of elliptic equations with a weighted gradient term}. Anal. PDE {\bf 9} (2016), 1671-1692.

\bibitem{NgV} Nguyen P. T., V\'{e}ron L., \textit{Boundary singularities of solutions to elliptic viscous Hamilton-Jacobi equations}. J. Funct. Anal. {\bf 263} (2012), 1487-1538.

\bibitem{PoQuSo} Polacik P., Quittner P., Souplet P.,  \textit{Singularity and decay
estimates in superlinear problems via Liouville-type theorems, Part.I:
Elliptic equations and systems}. Duke Math. J. {\bf 139} (2007), 555-579.

\bibitem{RiVe} Richard Y., V\'{e}ron L., \textit{Isotropic singularities of solutions of nonlinear elliptic inequalities}, Ann. Inst. H. Poincar\'e, {\bf 6} (1989), 37-72.

\bibitem{SeZo} Serrin J., Zou H., \textit{Existence and non-existence results for
ground states of quasilinear elliptic equations}. Arch. Ration. Mech. Anal. {\bf 121}
(1992), 101-130.

\bibitem{SeZo2} Serrin J., Zou H., \textit{Classification of positive solutions of quasilinear elliptic equations}. Topol. Methods Nonlinear Anal. {\bf 3}
(1994), 1-26.

\bibitem{Vo} Voirol F. X., {\it \'Etude de quelques \'equations elliptiques fortement non-lin\'eaires}. Th\`{e}se de Doctorat, Universit\'{e} de Metz (1994).


\end{thebibliography}
\end{document}